%% file: uniformization.tex
\documentclass{amsart}
\usepackage{amssymb}
\usepackage{eucal}
\usepackage{amsfonts}
\usepackage{epsfig}
\usepackage{color}
\usepackage[frame,ps,dvips,matrix,arrow,curve,rotate]{xy}
\let\oldlabel\label

\let\oldcite\cite

\let\oldmarginpar\marginpar
\renewcommand\marginpar[1]{\oldmarginpar{\scriptsize{#1}}}

\newtheorem{thm}{Theorem}[section]
\newtheorem{cor}[thm]{Corollary}
\newtheorem{lem}[thm]{Lemma}
\newtheorem{prop}[thm]{Proposition}
\theoremstyle{definition}
\newtheorem{defn}[thm]{Definition}
\theoremstyle{remark}
\newtheorem{rem}[thm]{Remark}
\numberwithin{equation}{section}
\theoremstyle{remark}

\newcommand{\X}{\mathcal{X}}
\newcommand{\Y}{\mathcal{Y}}

\newcommand{\PP}{\mathcal{P}}
\newcommand{\FF}{\mathcal{F}}
\newcommand{\U}{\mathcal{U}}
\newcommand{\V}{\mathcal{V}}
\newcommand{\W}{\mathcal{W}}

\newcommand{\D}{\mathcal{D}}
\newcommand{\C}{\mathcal{C}}
\newcommand{\Cm}{\C_{mod}}

\newcommand{\Xm}{\X_{mod}}

\newcommand{\mfG}{\mathfrak{G}}

\let\:=\colon

\newcommand{\Zn}{\mathbb{Z}_n}
\newcommand{\Zm}{\mathbb{Z}_m}
\newcommand{\Zd}{\mathbb{Z}_d}

\newcommand{\Ob}{\operatorname{Ob}}
\newcommand{\Mor}{\operatorname{Mor}}
\newcommand{\Out}{\operatorname{Out}}
\newcommand{\Aut}{\operatorname{Aut}}
\newcommand{\Hom}{\operatorname{Hom}}

\newcommand{\PSL}{\operatorname{PSL}}
\newcommand{\Pic}{\operatorname{Pic}}
\newcommand{\PGL}{\operatorname{PGL}}
\newcommand{\GL}{\operatorname{GL}}

\newcommand{\Ho}{\operatorname{Ho}}

\def\smashedlongrightarrow{\setbox0=\hbox{$\longrightarrow$}\ht0=1pt\box0}
\def\risom{\buildrel\sim\over{\smashedlongrightarrow}}
\def\smashedst{\setbox0=\hbox{$\rightrightarrows$}\ht0=4pt\box0}
\newcommand{\sst}[1]{\stackrel{#1}{\smashedst}}

\begin{document}

\title[Uniformization of Deligne-Mumford curves]
      {Uniformization of Deligne-Mumford curves}%

\author{Kai Behrend \and Behrang Noohi}%


\begin{abstract}
We compute the fundamental groups of non-singular analytic
Deligne-Mumford  curves, classify the simply connected ones, and
classify analytic Deligne-Mumford curves by their uniformization
type. As a result,
we find an explicit presentation of an arbitrary
Deligne-Mumford curve as a quotient stack.
Along the way, we  compute the automorphism 2-groups of
weighted projective stacks $\mathcal{P}(n_1,n_2,\cdots,n_r)$.  We also discuss
connections with the theory of $F$-groups, $2$-groups, and Bass-Serre theory
of graphs of groups.
\end{abstract}

\maketitle

\section{Introduction}

The aim of this paper is to extend the uniformization of Riemann
surfaces to  Deligne-Mumford analytic curves.

In the case of Riemann surfaces, the uniformization problem is
well-understood. There are three simply connected 1-dimensional
analytic domains, namely $\mathbb{H}$, $\mathbb{C}$ and
$\mathbb{P}^1$. Accordingly, there are three types of analytic
curves: {\em hyperbolic}, {\em Euclidian}, and {\em spherical}.
The only spherical curve is the Riemann sphere $\mathbb{P}^1$. The
Euclidian ones are the ones whose universal cover is $\mathbb{C}$.
These are precisely $\mathbb{C}$, $\mathbb{C}^*$, and elliptic
curves. Everything else has universal cover $\mathbb{H}$, so is
hyperbolic.

Every Riemann surface arises as a quotient of one of these  simply
connected domains by a discrete subgroup of the group of its
analytic automorphisms. Such discrete subgroups, being the
fundamental groups of the corresponding Riemann surfaces, have a
particularly simple form, and they act without fixed points on the
corresponding symmetric domain.

In passage to Deligne-Mumford curves, we face some complications,
as there are more simply connected Deligne-Mumford curves which
were not present in the classical case. For instance, every
weighted projective line $\PP(m,n)$ is simply connected. The good
news is that these are all.

\begin{thm}
    The simply connected Deligne-Mumford
curves are precisely ${\mathbb H}$, ${\mathbb C}$ and ${\mathcal
P}(m,n)$, for arbitrary $n,m\geq1$.
\end{thm}

So, the trichotomy hyperbolic/Euclidian/spherical continues to
hold in the more general Deligne-Mumford setting: every
Deligne-Mumford curve has a universal cover which is either
${\mathbb H}$, ${\mathbb C}$ or ${\mathcal P}(m,n)$, for some
$n,m\geq1$.

The homotopy groups of a Deligne-Mumford curve are rather easy to
compute. It is also easy to determine the uniformization type of a
Deligne-Mumford curve by just ``looking at it''! (Propositions
\ref{P:type}, \ref{P:hEs}.)

The next task is to represent a Deligne-Mumford curve as a
quotient of its universal cover. This is pretty easy in the {\em
hyperbolic} and {\em Euclidian} cases (Proposition \ref{P:hE}).
The spherical case, however, poses some challenge. The key
difference is that, in the spherical case the universal cover is
itself a stack, so its automorphisms from a {\em 2-group}. This
means, to be able to talk about deck transformations, we have to
incorporate some 2-group theory. As a result the classification of
spherical Deligne-Mumford curves is less trivial. The main inputs
are: 1)  computation of the automorphism 2-groups of weighted
projective lines, 2) homotopy classification of maps from discrete
groups into these 2-groups. After accomplishing this, we arrive at
the following classification results.

\begin{thm}
  Let $m\neq n$ be given positive integer numbers, and let
  $d=\gcd(m,n)$. Then, the isomorphism classes of Deligne-Mumford
  curves $\X$ with universal cover $\PP(m,n)$ are in natural bijection with
  isomorphism
  classes of pairs $(K,\chi)$, where $K$ is a finite group containing
  $\mu_d$ as a central subgroup and
  $\chi \: K/\mu_d \to \mathbb{C}^*$ is a character.
\end{thm}

\begin{thm}
 All Deligne-Mumford curves whose universal cover is ${\mathcal P}(d,d)$
 are given by triples $(\iota,E,\rho)$, where
  \begin{itemize}
   \item[$\mathbf{i.}$]
    $\iota:{\mathbb C}^\ast\hookrightarrow E$ is a
     central extension of algebraic groups,
   \item[$\mathbf{ii.}$] $\rho:E\to GL_2$ is a 2-dimensional representation,
  \end{itemize}
 such that the following two conditions are satisfied:

  \begin{itemize}
   \item[$\mathbf{a.}$]  $\rho\circ\iota:{\mathbb C}^\ast\to GL_2$
      is the diagonal $(d,d)$-th power map,
   \item[$\mathbf{b.}$] $E/{\mathbb C}^\ast$ is finite.
  \end{itemize}
\end{thm}

 The above results are constructive. We show in Section \ref{S:spherical}
 how given data as in above theorems one can produce the corresponding
Deligne-Mumford curve as a quotient stack. In particular, this
gives a canonical presentation of a spherical Deligne-Mumford
curve as a quotient stack. (We have such a presentation in the
hyperbolic and Euclidian case too.)

The paper is organized as follow. In Section \ref{S:NT} we fix
some notations. In Section \ref{S:Analytic} we overview the theory
of analytic stacks. Section \ref{S:Pre} serves as a toolbox
throughout the paper.  In it we collect some easy, and perhaps
well-known, facts which will be  invoked in the course of proofs
or constructions in the ensuing sections. This is to avoid
cluttering up the arguments with  trivialities.

Section \ref{S:UDMO} discusses the uniformization problem in the
orbifold case. We compute (using van Kampen theorem) the
fundamental groups of all orbifold curves, classify the simply
connected ones, and classify the orbifold curves by their
uniformization type. In Section \ref{S:UDMG} we extend this to all
smooth Deligne-Mumford curves. In this case we have a new family
of simply connected stacks which should be taken into
consideration. These are the weighted projective lines. In
Subsection \ref{SS:graphs} of Section \ref{S:UDMG} we  give a  way
of computing the fundamental group of a Deligne-Mumford curve
avoiding the theorem of van Kampen. The point is that, an (open)
Deligne-Mumford curve contains a graph of groups \cite{Serre} as a
1-skeleton that is a deformation retract. In particular the
formulas in \cite{Serre} give an alternative way of computing the
fundamental groups.

Section \ref{S:spherical} augments Section \ref{S:UDMG}, or more
precisely,  Proposition \ref{P:spherical2} of Section
\ref{S:UDMG}, by giving a more conceptual explanation of the
theorem and describing the classification of spherical
Deligne-Mumford curves (the ones whose universal cover is a
weighted projective line).

\vspace{0.1in}

\noindent{\bf Acknowledgments.} The second author would like to
thank University of British Columbia, Newton Institute, and
C.I.R.M-Luminy, where the research for this work done, as well as
B. Toen and A. Vistoli for helpful discussions.

\tableofcontents

\section{Notations and terminology}{\label{S:NT}}

We use the terms {\em Riemann surface} and {\em analytic curve}
interchangeably. All Riemann surfaces and Deligne-Mumford curves
considered in this paper are assumed to be of finite type (i.e.
compact with a finite number of discs, possibly of zero radius,
removed).

The symbols $\mathbb{H}$, $\mathbb{C}$ and $\mathbb{P}^1$ stand
for the upper half-plane, the complex plane and the Riemann
sphere, respectively. We denote the complex unit disc by
$\mathbb{D}$ and its stacky quotient under the rotation action of
$\Zn$ by $\D_n$. We denote the punctured disc $\mathbb{D}-\{0\}$
by $\D_{\infty}$ or $\mathbb{D}^*$. Similarly, we use the
notations $\C_n$ for $[\mathbb{C}/\Zn]$, and $\C_{\infty}$ or
$\mathbb{C}^*$ for $\mathbb{C}-\{0\}$. The complex plane with two
orbifold points of order 2, say at $0$ and $1$, is denoted by
${\mathcal C}_{2,2}$.

All stacks are assumed to be smooth and separated. By an {\em
orbifold} we mean a (smooth) Deligne-Mumford (topological,
analytic or  algebraic) stack that is generically a (topological,
analytic, algebraic) space. (Some authors call this a {\em reduced
orbifold}.) By an {\em orbifold point} of a stack (which is not
necessarily an orbifold) we mean a point at which the inertia
group jumps (e.g. origin in $\D_n$,  $2 \leq n < \infty$).

We denote the coarse moduli space (sometime called the
``underlying space") of a (topological, analytic or algebraic)
stack $\X$, $\C$ etc. by $\Xm$, $\Cm$ etc..

We do not want to worry about 2-isomorphisms so 2-isomorphic
morphisms are declared equal, except in Sections \ref{S:wtd} and
\ref{S:spherical}. When talking about stacks, we use the terms
{\em equivalence} and {\em isomorphism} interchangeably.

We say that a topological stack is {\em uniformizable}, if its
universal cover is an honest topological space. We point out that
the universal cover of a topological stack always exists as a
topological stack. This is proven in \cite{homotopy1}.

For a subgroup $H$ of a group $G$, we denote its centralizer by
$Z_G(H)$. We denote the center of $G$ by $Z(G)$.

We denote the general linear group over $\mathbb{C}$ by
$\GL_n(\mathbb{C})$, or more briefly, $\GL_n$. Similarly $\PGL_n$
stands for the projective general linear group over $\mathbb{C}$.

Sometime we use the term `map' instead of `morphism'. For
instance, instead of saying a `morphism $\Gamma \to \mathfrak{G}$
of 2-groups' we say a `map $\Gamma \to \mathfrak{G}$ of 2-groups'.

\section{Overview of topological and analytic stacks} {\label{S:Analytic}}

We recall the basic definitions/facts about topological and
analytic stacks. For more details the reader can consult
\cite{homotopy1}.

\subsection{Topological stacks}
Let $\mathbf{Top}$ be the category of topological spaces. We
endow $\mathbf{Top}$ with the usual Grothendieck topology (covers
are simply topological open covers); so we can talk about (the
2-category of) stacks over $\mathbf{Top}$. By Yoneda's lemma, the
category of topological spaces embeds as a full subcategory of the
2-category of stacks.

We say a morphism $f \: \Y \to \X$ of stacks over  $\mathbf{Top}$
is {\em representable}, if for any map $X \to \X$ from a
topological space $X$ to $\X$, the fiber product
$Y=X\times_{\X}\Y$ is equivalent to a topological space. Let
$\mathbf{P}$ be a property of maps of topological spaces that is
stable under base change; for example $\mathbf{P}$=local
homeomorphism, $\mathbf{P}$=open, $\mathbf{P}$=finite fibers,
$\mathbf{P}$=discrete fibers, or $\mathbf{P}$=covering map. We say
a representable morphism $f \: \Y \to \X$ of stacks over
$\mathbf{Top}$ satisfies $\mathbf{P}$ if for any map $X \to \X$
from a topological space $X$ to $\X$, the base extension $Y \to X$
of $f$ satisfies $\mathbf{P}$. A {\em pre-Deligne-Mumford
topological stack} is  a stack for which there exists an
epimorphism $p \: X \to \X$ from a topological space $X$, such
that $p$ is representable by local homeomorphisms.

Sometimes we need to work with properties $\mathbf{P}$ that are
only invariant under base extension via {\em local homeomorphism}.
For example $\mathbf{P}$=closed. In this case, we can still define
the property $\mathbf{P}$ for a morphisms $f \: \Y \to \X$ of
stacks, provided $\X$ is a pre-Deligne-Mumford stack. More
precisely, we say that $f \: \Y \to \X$ satisfies $\mathbf{P}$, if
for any local homeomorphism $X \to \X$ from a topological space
$X$ to $\X$, the base extension $Y \to X$ of $f$ satisfies
$\mathbf{P}$.

 \begin{defn}{\label{D:topologicalDM}}
   A pre-Deligne-Mumford topological stack $\X$ is called a {\bf Deligne-Mumford
   topological stack} if the diagonal $\X \to \X\times\X$ is representable by
   closed maps with discrete finite fibers.
 \end{defn}

 \begin{rem}
   In \cite{homotopy1} the finiteness condition in the definition
   of a Deligne-Mumford topological stack is replaced by a weaker
   condition which allows infinite inertia groups.
 \end{rem}

 With some extra care, the natural constructions of homotopy theory
 of topological spaces can be extended to Deligne-Mumford topological stacks.
 For instance, one can define a reasonable notion of homotopy between maps,
 and this allows us to define the $n$-th homotopy group
 of a pointed topological stack $(\X,x)$ as pointed homotopy classes of maps
 $[(S^n,*),(\X,x)]$.

 In particular, we can talk about the fundamental group
 $\pi_1(\X,x)$ of a topological stack. When $\X$ is locally well-behaved (i.e.
 it is locally path-connected and semilocally 1-connected), the fundamental
 group classifies the connected covering spaces of $\X$
 (that is, there is a natural
 bijection between the equivalence classes  of pointed connected covering spaces
 of $\X$ and subgroups of $\pi_1(\X,x)$; in particular, $\X$ has a universal
 cover $\tilde{\X}$ that is unique up to equivalence, and $\X$ is the quotient
 of $\tilde{\X}$ by a natural action of $\pi_1(\X,x)$).
 Recall from \cite{homotopy1} that a covering space of $\X$
 is defined to be a representable map  $\Y \to \X$ of stacks
 that is $\mathbf{P}$, where $\mathbf{P}$=covering map (see the second
 paragraph of this section). Note that we have not assumed $\Y$ to be connected.
 For more  details on the covering theory of stacks the reader may
 consult \cite{homotopy1}.

 To any topological stack $\X$ one can associate an honest topological space
 $\Xm$, called the  {\em coarse moduli space}, or more informally the {\em
 underlying space}, of $\X$. The coarse moduli space $\Xm$ is in some sense the
 best approximation of $\X$ by a topological space. For instance, when $\X$
 is a quotient stack of a group action, then $\Xm$ is just the usual coarse
 quotient. There is a natural continuous map $\X \to \Xm$, called the
 {\em moduli map}. Usually we visualize a topological stack by drawing its
 coarse moduli space and describing the stacky structure by specifying
 the inertia groups of the points; of course this naive picture
 is not enough to recover
 the stack. For example, given a loop on $\X$, we draw it as a loop on $\Xm$,
 but it should be kept in mind that different loops in $\X$ may give
 rise to the same drawing. For example, typically we have non-trivial
 loops in $\X$ that become constant loops in $\Xm$.

\subsection{Complex analytic stacks}
Since we are only interested in smooth analytic stacks, we give a
simplistic definition of (smooth) Deligne-Mumford analytic stacks
which makes life easier for those who do not feel comfortable with
analytic spaces (see Remark \ref{R:analytic} below).

Let $\mathbf{Comp}$ be the category of complex manifolds, endowed
with the usual Grothendieck topology (where covers are simply
topological open covers). As in the case of topological stacks, we
can construct the 2-category of stacks over $\mathbf{Comp}$, and,
by Yoneda,  this category contains the category of complex
manifolds as a full subcategory.

We say a morphism $f \: \Y \to \X$ of stacks over  $\mathbf{Comp}$
is {\em representable by local homeomorphisms}, if for any map $X
\to \X$ from a complex manifold $X$ to $\X$, the fiber product
$Y=X\times_{\X}\Y$ is equivalent to a complex manifold, and the
map $Y \to X$ is a local homeomorphism. A stack $\X$ over
$\mathbf{Comp}$ is called a (smooth) {\em pre-Deligne-Mumford
analytic stack} if there exists an epimorphism $p \: X \to \X$
from a complex manifold $X$ such that $p$ is representable by
local homeomorphisms.
 We say a morphism $f \: \Y \to \X$
of pre-Deligne-Mumford analytic stacks is {\em representable}, if
for any map $X \to \X$ from a complex manifold $X$ to $\X$ that is
representable by local homeomorphisms, the fiber product
$Y=X\times_{\X}\Y$ is equivalent to a complex manifold. Let
$\mathbf{P}$ be a property of morphisms of complex manifolds that
is invariant under base change along local homeomorphisms. For
example we can take $\mathbf{P}$=closed,   $\mathbf{P}$=finite
fibers, $\mathbf{P}$=discrete fibers or $\mathbf{P}$=covering
space. Then, we say a representable morphism $f \: \Y \to \X$ of
pre-Deligne-Mumford analytic stacks is $\mathbf{P}$, if for any
map $X \to \X$ from a complex manifold $X$ to $\X$ that is
representable by local homeomorphisms, the base extension $Y \to
X$ is $\mathbf{P}$.

 \begin{defn}{\label{D:analyticDM}}
   A pre-Deligne-Mumford analytic stack $\X$ is called a {\bf Deligne-Mumford
   analytic stack} if the diagonal $\X \to \X\times\X$ is representable by
   closed maps with finite fibers.
 \end{defn}

 \begin{rem}{\label{R:analytic}}
  In \cite{homotopy1} the definition of a Deligne-Mumford analytic stack is
   different from the one given above in two ways. Firstly, here we have added
   the finiteness condition  on the diagonal  to make the exposition simpler.
   Secondly, in
   \cite{homotopy1} the base Grothendieck site is taken to be the category
   $\mathbf{Analytic}$ of all analytic spaces,instead of the smaller
   $\mathbf{Comp}$. However, apart from the finiteness conditions on
   the diagonal, the two theories are ``the same" if
  we restrict out attention to {\em smooth} analytic stack. In other words,
  the
  2-category of analytic stacks defined above is equivalent to the sub 2-category
  of the 2-category of analytic stack of \cite{homotopy1} consisting of
  smooth analytic stacks. In more technical terms, if $\X$ is a smooth analytic
  stack in the sense of \cite{homotopy1}, one can always recover $\X$ from
  its restriction to $\mathbf{Comp}$ using a homotopy left Kan extension.
  More casually speaking, this means that, for a smooth $\X$,
  if we know the maps from
  {\em smooth} domains into $\X$, then we know the  maps from an arbitrary
  domain into $\X$.
 \end{rem}

 We quote the following result from \cite{homotopy1}.

  \begin{prop}{\label{P:localquotient}}
    Let $\X$ be a Deligne-Mumford topological (resp., analytic) stack. Then
    there is a covering $\{\U_i\}$ of $\X$ by open substacks such that
    each $\U_i$ is a quotient stack $[X/G]$, where $X$ is a topological space
    (resp., complex manifold), and $G$ a finite group acting continuously
    (resp., analytically) on $X$.
  \end{prop}

 One can define the coarse moduli space $\Xm$ of an analytic Deligne-Mumford
 stack $\X$. It is an analytic space, but it may not in general be a smooth complex
 manifold. When $\X$ is 1-dimensional, however, which is the only case  we are interested
 in anyway, the coarse moduli space $\Xm$ is
 an analytic curve (i.e., a Riemann surface).

\subsection{Comparing algebraic, analytic and topological stacks}
  In any of the algebraic, analytic or topological settings,
  there is an alternative description of Deligne-Mumford stacks
  using groupoids. For instance, in the algebraic setting, a
  Deligne-Mumford stack $\X$ can be described by (the Morita equivalence
  class of) an \'{e}tale groupoid $X_1\sst{}X_0$, where the diagonal
  $X_1\sst{}X_0\times X_0$ is assumed to be a finite
  morphism.
  We have similar descriptions in the analytic and topological settings.
  This enables us to define natural functors
    $$\mathbf{AlgDM \to AnDM \to TopDM}.$$
  For instance, let $\X$ be a smooth algebraic Deligne-Mumford stack of finite
  type over complex numbers, with
  $X_1\sst{}X_0\times X_0$ an \'{e}tale groupoid representing it. Then we
  define $\X^{an}$ to be the quotient of the groupoid
  $X_1^{an}\sst{}X_0^{an}\times X_0^{an}$.  

  We can do the same thing to get a topological stack $\X^{top}$ from an
  analytic (or algebraic) stack $\X$. The homotopy groups of an
  analytic (or algebraic) stack $\X$ are defined to be those of $\X^{top}$.

  The proof that the above functors are well-defined can be found in
  \cite{homotopy1}, where the reader can also find the proof of the fact that
  these functors respect the coarse moduli space construction.

\section{Some useful facts} {\label{S:Pre}}

\subsection{Basic facts about fundamental groups}{\label{SS:BFFG}}

In \cite{Noohi} it is shown that for any pointed algebraic stack
$(\X,x)$ there is a natural group homomorphism $I_x \to
\pi_1(\X,x)$, where the group $I_x$ is   the inertia group of the
residue gerbe at $x$.\footnote{In \oldcite{Noohi} the groups $I_x$
are named {\em hidden fundamental groups} and are denoted by
$\pi_1^h(\X,x)$.}
 Similar statement for topological (hence
also analytic) stacks is proved in \cite{homotopy1}. We quote the
following results from {\it loc.\,cit.}:

\begin{lem}{\label{L:hidden}}
Let $\X$  be a Deligne-Mumford topological stack, and let $(\Y,y)
\to (\X,x)$ be a pointed  covering space. Then we have a cartesian
diagram
      $$\xymatrix{  I_y \ar[r] \ar[d] &   \pi_1(\Y,y) \ar[d] \\
                    I_x \ar[r]        &   \pi_1(\X,x)}$$
\noindent In particular, the kernel of $I_x \to \pi_1(\X,x)$ is
invariant under covering maps. The same statement is true for
analytic and algebraic stacks (where in the algebraic situation
$\pi_1(\X,x)$ refers to algebraic fundamental group, and covering
maps are understood to be finite \'{e}tale maps).
\end{lem}

\begin{thm}{\label{T:uniformizable}}
Suppose $\X$ is a Deligne-Mumford topological stack.
 \begin{itemize}
   \item[$\mathbf{i.}$]
     Let $\Y \to \X$ be a covering space,
     and let $H \subseteq \pi_1\X$ be the corresponding subgroup (for a choice
     of base points). Then, $\Y$ is an honest topological space if and only
     if for any point $x$, the map $I_x \to \pi_1(\X,x)$ is
     injective and its image does not intersect any conjugate of $H$.

   \item[$\mathbf{ii.}$]  $\X$ is
       uniformizable (i.e.  its universal cover is an honest topological space)
       if and only if the maps   $I_x \to \pi_1(\X,x)$ are
       injective for every  point $x$ in $X$. If, furthermore, the maps
       $I_x \to \pi_1(\X,x)$ remain injective after passing to profinite
       completion, then $\X$ is the quotient stack of a finite group
       acting on a topological space.
 \end{itemize}

 These statements remain valid for
analytic and algebraic stacks (where in the algebraic situation
$\pi_1(\X,x)$ refers to the algebraic fundamental group and
uniformizable means that $\X$ is the quotient stack of the action
of a finite group on  an algebraic space).
\end{thm}

\begin{rem}
In this paper the only algebraic spaces that will come up are
smooth and 1-dimensional, in which case they are automatically
schemes.  So in the algebraic version of the previous theorem we
will actually get smooth algebraic curves covering $\X$.
\end{rem}

The following theorem is proved in \cite{homotopy1}.

\begin{thm}{\label{T:RET}}
  Let $\X$ be an algebraic stack of finite type over $\mathbb{C}$.
  Then the algebraic fundamental group of $\X$ is isomorphic to the
  profinite completion of the topological fundamental group of the underlying
  topological stack $\X^{top}$.
\end{thm}

The van Kampen theorem remains true for topological stacks. The
proof will appear elsewhere.

\subsection{Stacky discs and their fundamental groups}{\label{SS:SD}}
In this part we look at (smooth) Deligne-Mumford analytic stacks
whose coarse moduli space is isomorphic to  $\mathbb{D}$ or
$\mathbb{C}$. More specifically, we are interested in ones that
are gerbes away from the origin.

The following two propositions explains why we would like to study
stacky discs.

 \begin{prop}{\label{P:Dn}}
   Let $\U$ be a Deligne-Mumford analytic orbifold whose coarse
   moduli space is $\mathbb{D}$ (resp. $\mathbb{C}$). Assume
   $\U$ has at most one orbifold point. Then $\U$ is isomorphic to
   $\D_n$ (resp. $\C_n$) for some $n \geq 1$.
   (See Section \ref{S:NT} for notation.)
 \end{prop}

  \begin{proof}
    We consider the case where the coarse moduli space of $\U$
    is $\mathbb{D}$; the proof goes through similarly in the
    case where it is $\mathbb{C}$. We will assume that $\U$
    is not a space (so it has a unique orbifold point).

    Let $x$ be the orbifold point of $\U$.
        By Proposition \ref{P:localquotient}, there is an open substack
    $\U' \subset \U$ around $x$ that is a quotient stack by a finite group
    action. After shrinking $\U'$, we may assume that its underlying space
    is a disc. Now, by Theorem \ref{T:uniformizable},
    $I_x\to \pi_1(\U',x)$ is injective. But $\pi_1(\U',x)=\pi_1(\U,x)$.
     So $I_x\to \pi_1(\U,x)$
    is injective. Hence, by  Theorem \ref{T:uniformizable}, $\U$
    is uniformizable. The universal cover of $\U$ is necessarily either
    $\mathbb{D}$ or $\mathbb{C}$. But $\mathbb{C}$ is ruled out, since
    composing the map $\mathbb{C} \to \U$ with $\U \to \U_{mod}\cong\mathbb{D}$
    would give a surjective map $\mathbb{C} \to \mathbb{D}$ which is impossible.
    So we have $\U \cong [\mathbb{D}/G]$, where $G$ is a discrete group
    acting properly
    discontinuously on $\mathbb{D}$. Since $\U$ has a unique orbifold
    point, this action has a unique orbit $\mathcal{O}$
    whose elements have non trivial
    stabilizers. If we remove this orbit from $\mathbb{D}$, we find a
    covering space
    of the punctured disc $\mathbb{D}^*$. But the only possible
    covering spaces of $\mathbb{D}^*$ are isomorphic to either
    $\mathbb{D}$ or $\mathbb{D}^*$. This implies that the orbit
    $\mathcal{O}$ is in fact just a single point, which after some
    change of coordinates  can be taken to be $0 \in \mathbb{D}$. So,
    $0$ is the unique fixed point of the action of $G$.
    Since $\U$ is a Deligne-Mumford stack, $G$ must be finite.
    The only finite groups of automorphisms of $\mathbb{D}$
    which fix the origin are the cyclic rotation groups. So $G \cong
    \mathbb{Z}_n$,
    for some $n \geq 2$, and it acts by rotations around origin. Therefore,
    $\U\cong\D_n$.
 \end{proof}

 This proposition implies that an analytic orbifold curve $\C$ is
  uniquely defined by
  giving a Riemann surface $C$ (to be $\Cm$), a finite collection $A$
  of points on $C$ (the orbifold points),
  and an integer $n_x \geq 2$ for each $x \in A$ (the orbifold degree at $x$).

 \begin{prop}{\label{P:gerbe}}
  Let $\X$ be a Deligne-Mumford  analytic stack. Then $\X$ is an $H$-gerbe
  over an orbifold $\Y$, for some finite group $H$. Furthermore,
  the construction of $\Y$ is functorial in $\X$, in the sense that,
  an arbitrary map $\X' \to \X$ of Deligne-Mumford  analytic stacks
  induces a natural map $\Y' \to \Y$ of the corresponding orbifolds,
   making the following diagram commute:
  $$\xymatrix@=12pt@M=10pt{
   \X' \ar[r]\ar[d]  &  \X \ar[d]  \\
            \Y' \ar[r]        &   \Y    }$$
   The same statement
  is true for algebraic and topological Deligne-Mumford stacks. (Here $H$-gerbe
  means, locally (in the \'{e}tale topology) our gerbe is equivalent to
  the classifying stack of $H$.)
 \end{prop}

  \begin{proof}
    We only give the proof
     for the analytic  stacks because that is
    what we need. The exact same proof works for topological stacks, and the
    same  idea can be made to work in the algebraic case too.

    Let $R \sst{} X$ be an \'{e}tale groupoid representing $\X$,
    and let $S \subseteq R$ be its stabilizer group. We claim that there
    is unique maximal locally constant subgroup $T$  of $S$ (as groups over
    $X$), and that $T$ is normal in $R$. Having shown this, it is enough to
    consider the groupoid $R' \sst{} X$, where $R'=R/T$ is the quotient
    groupoid, and define $\Y$ to be $[X/R']$.

    Since the statement is local we can use Proposition
    \ref{P:localquotient} and assume that $\X=[X/G]$, where $G$ is a finite
    group acting on a complex manifold. Let $H\subseteq G$ be the subgroup
    of elements that act trivially on $X$. It is easy to see
    that $H$ is normal in $G$ and that, as a subgroupoid of the stabilizer of
    the groupoid $G\times X \sst{} X$, it is the largest
    locally constant subgroup.
    This proves the claim. The proof of functoriality is easy and is also
    left to the reader.
  \end{proof}

 Now, suppose $\U$ is a Deligne-Mumford analytic  stack whose coarse moduli
 space is $\mathbb{D}$,
and that, away from the origin,  is an $H$-gerbe over
$\mathbb{D}^*$ for some finite group $H$. Let $G$ be the inertia
group of $\U$ at the origin. By Proposition \ref{P:gerbe} and
Proposition \ref{P:Dn},
 $\U$ is an $H$-gerbe over some $\D_n$, $n \geq 1$.

We have a short exact sequence
                  $$ 1 \to H \to G \to \Zn \to 1, $$
\noindent and this short exact sequence uniquely characterizes
$\U$ up to isomorphism (Proposition \ref{P:characterize}). Here is
how we can reconstruct $\U$ from this short exact sequence: let
$G$ act on $\mathbb{D}$ via its quotient $\Zn$ (by rotations). The
quotient of this action is isomorphic to $\U$. In particular, $\U$
is uniformizable, and its fundamental group is isomorphic to $G$.
More generally, we have the following.

\begin{prop}{\label{P:characterize}}
  Let $\Y$ be a topological stack such that $\pi_n\Y$ is trivial
  for $n\geq 2$. Let $H$ be a discrete group. Let $\X$ be an $H$-gerbe
  over $\Y$. Then, there is a short exact sequence
    $$ 1 \to H \to \pi_1\X \to \pi_1\Y \to 1 $$
  and this short exact sequence characterizes $\X$ up to isomorphism.
  Furthermore, $\X$ and $\Y$ have the same universal cover.
  In particular, $\X$ is uniformizable if and only if $\Y$ is.
  In this case, we have $\X\cong[\tilde{\Y}/\pi_1\X]$, where
  $\pi_1\X$ acts on $\tilde{\Y}$ via the map $\pi_1\X \to \pi_1\Y$).
\end{prop}

\begin{proof}
  The short exact sequence is simply the fiber homotopy exact sequence.
  To prove that  it characterizes $\X$, it is enough to
  show that $\X$ and $\Y$ have the same universal cover (because then
  $\X$ will have to be the quotient stack $[\tilde{\Y}/\pi_1\X]$).

  To prove that $\X$ and $\Y$ have the same universal cover let us consider
  the following commutative diagram:
  $$\xymatrix@=12pt@M=10pt{
     \tilde{\X} \ar[r]^q\ar[d]  &  \tilde{\Y} \ar[d] \\
           \X  \ar[r]_p        &    \Y   }$$
  By the 2-out-of-3 property of finite \'{e}tale maps, it follows
  that $q$ is   finite \'{e}tale (not, a priori, representable).
  Let $\FF$ be the fiber of $q$. Since $\FF$ is  finite \'{e}tale over a point,
  it can only have non-trivial $\pi_0$ and $\pi_1$. On the other hand, since
  all homotopies of $\tilde{\Y}$ vanish, the fiber homotopy exact sequence
  implies that the map $\FF \to \tilde{\X}$ induces
  an isomorphism on all $\pi_i$, $i \geq  0$. Since $\tilde{\X}$ is connected
  and simply connected, this implies that $\pi_i \FF$ is trivial for $i \geq  0$.
  That is, $\FF$ is just a point. Therefore, $q$ is an isomorphism.
\end{proof}

 Observe that the short exact sequence of the proposition
 gives rise to an element
in $\Hom(\pi_1\Y,\Out{H})$, the {\em band} of the gerbe. When
$\Y=U$ is an open (non-orbifold) analytic curve, the second
cohomology group (with coefficients in an abelian sheaf) is
trivial, so  an $H$-gerbe over $U$ is characterized  by its
band.\footnote {Having  fixed an isomorphism $\pi_1U \cong
\mathbb{Z}^{*m}$, this amounts to specifying $m$ elements in
$\Out{H}$. We will use these elements to construct the graph of
groups associated to an open Deligne-Mumford curve (see Section
\ref{SS:graphs})}

\vspace{0.2in}
\noindent{\bf Gerbes with trivial band.} Now let us look more
closely at  gerbes over $\D_n$ whose band is trivial. The
following sets are in natural bijection:
\begin{itemize}
  \item $H$-gerbes over $\D_n$ whose band is trivial (with a
          chosen trivialization), up to isomorphism.
  \item Extensions of $\Zn$ by $H$ for which the centralizer of $H$ surjects
      onto $\Zn$, up to an isomorphism inducing identity on $H$ and $\Zn$.
  \item Central extensions of $\Zn$ by $Z(H)$, up to an isomorphism inducing
        identity on $Z(H)$ and $\Zn$.
  \item $Z(H)/nZ(H)$.
\end{itemize}
Let us be more explicit about  these correspondences because we
will use them later for calculating some fundamental groups.
Consider an $H$-gerbe over $\D_n$  with trivial band. The
corresponding group extension
                    $$ 1 \to H \to G \to \Zn \to 1 $$
\noindent has the property that for any element in $\Zn$ the
corresponding outer action on $H$ is trivial (this means the band
is trivial). Therefore, $G=Z_G(H)H$. This gives rise to the
following central extension
                    $$ 1 \to Z(H) \to Z_G(H) \to \Zn \to 1.$$
\noindent We can recover the original extension by pushing out
along $Z(H) \hookrightarrow H$. Let us now elaborate on the
bijection between the extensions of $\Zn$ by $H$ and
$Z(H)/nZ(H)$. Given an element $a \in Z(H)$, define
                    $$G_a:=H \times \mathbb{Z} /(a,-n).$$
We have an exact sequence
                    $$ 1 \to H \to G_a \to \Zn \to 1,$$
\noindent which is, up to an isomorphism inducing identity on $H$
an $\Zn$, independent of the class of $a$ in $Z(H)/nZ(H)$.
Conversely, given an extension
                    $$ 1 \to H \to G \to \Zn \to 1$$
\noindent pick an element $g \in G$ which maps to $1 \in \Zn$. The
conjugation action of $g$  on $H$ is inner by assumption. So we
can find $h \in H$ which induces the same action. Let $a:=h^n$.
The conjugation action of $a$ on $H$ must then be trivial. That
is, $a$ lies in the center of $H$. The class of $a$ in
$Z(H)/nZ(H)$ is easily seen to be independent of the choice of
$g$.

\vspace{0.2in}
\noindent{\bf Extension of a gerbe from $\D_{\infty}$ to $\D_n$.}
It is an easy exercise to check that an $H$-gerbe over
$\D_{\infty}$ extends to an $H$-gerbe over $\D_n$ if and only if
the corresponding element in $\Out(H)$ is $n$-torsion. The
extensions of a given gerbe are acted upon simply transitively by
$H^2(\Zn,Z(H))$. Here $Z(H)$ is made into a   $\Zn$-module via the
action induced by the given $n$-torsion element in $\Out(H)$. We
have the corresponding cartesian diagram of group extensions:
{\footnotesize
$$\xymatrix{ 1 \ar[r] & H \ar[r]\ar[d]_{=} & G_0 \ar[d] \ar[r] &  \mathbb{Z} \ar [r] \ar[d] & 1 \\
             1 \ar[r] & H \ar[r]           & G  \ar[r]         &    \Zn  \ar[r]             & 1 }$$
} \noindent Note that $G_0 \cong H \rtimes \mathbb{Z}$. The
semidirect product is defined by choosing a lift of the given
element of  $\Out(H)$ to $\Aut(H)$. This isomorphism becomes
canonical once we choose a lift of  $1 \in
\mathbb{Z}=\pi_1{\D_{\infty}}$ to a loop in the gerbe.

\vspace{0.2in}
\noindent{\bf Differential forms on $\D_n$.} The sheaf of
differentials on $\D_n$ is a line bundle. Sections of this line
bundle are the  $\Zn$-invariant differential forms on
$\mathbb{D}$, necessarily of the form $z^{n-1}f(z^n)dz$.
Therefore, $z^{n-1}dz$ is a generator of the sheaf of
differentials on $\D_n$. Its order of vanishing is
$\frac{n-1}{n}$.

\subsection{Weighted projective lines, footballs, and drops}{\label{SS:FD}}
For any pair of integers $m,n \geq 1$ there is a Deligne-Mumford
curve $\PP(m,n)$, called the {\em weighted projective line of type
$(m,n)$}. It has at most two orbifold points, and its coarse
moduli space is $\mathbb{P}^1$. Here is the standard way of
constructing $\PP(m,n)$:\footnote{See Section \ref{S:GoverF} for
an alternative construction} consider the action of
$\mathbb{C}^{*}$ on $\mathbb{C}^2- \{0\}$ given by
$t.(x,y)=(t^mx,t^ny)$. The corresponding quotient stack is
$\PP(m,n)$. In particular, we have  a $\mathbb{C}^{*}$-fibration
over our weighted projective line whose total space is
$\mathbb{C}^2-\{0\}$. A fiber homotopy exact sequence argument
shows that $\PP(m,n)$ is simply connected, for all $m,n \geq 1$.
In fact, $\pi_i\PP(m,n)\cong\pi_i S^2$ for all $i$.

Obviously, $\PP(m,n)\cong\PP(n,m)$. A weighted projective line is
an orbifold if an only if $m$ and $n$ are relatively prime. We
call these {\em orbifold weighted projective lines}. Note that
$\PP(1,1)=\mathbb{P}^1$.    The weighted projective line
$\PP(m,n)$ is a $\Zd$-gerbe over the weighted projective line
     $\PP(\frac{m}{d},\frac{n}{d})$, where $d$ is the greatest common divisor.

For a pair of integers $m,n \geq 1$, a {\em football} of type
$(m,n)$ is an orbifold whose underlying space is $\mathbb{P}^1$
and has two orbifold points of order $m$ and $n$. We denote by
$\FF(m,n)$ the football whose underlying space is $\mathbb{P}^1$
and has orbifold points of order $m$ and $n$ at $0$ and $\infty$,
respectively. Any football of type $(m,n)$ is isomorphic to
$\FF(m,n)$. When exactly one of $m$ or $n$ is $1$, a football of
type $(m,n)$ is sometimes called a {\em drop}.

When $m$ and $n$ are relatively prime, we have a natural
isomorphism $\FF(m,n)\cong\PP(m,n)$. In other cases, $\FF(m,n)$
and $\PP(m,n)$ are non-isomorphic.

When $m$ and $n$ are strictly greater than $1$, the group of
automorphisms of $\FF(m,n)$ is $\mathbb{C}^{*}$, the action being
the obvious rotation action. When exactly one of $m$ and $n$ is
equal to $1$, say $m=1$, the group of automorphisms of $\FF(m,n)$
is naturally isomorphic with the group of affine transformations
of the complex plane $\mathbb{C} \subset \FF(1,n)$.

\section{Uniformization of Deligne-Mumford curves:
                          the case of an orbifold}{\label{S:UDMO}}

In this section we compute the fundamental groups of  orbifold
analytic curves, determine the simply connected ones, and classify
orbifold curves by their uniformization types. The results of this
section are quite elementary,  presumably well-known (see for
instance \cite{FK} Section IV.9; especially, Theorem IV.9.12).

 We will assume that our orbifold curves $\C$ are of finite type.
 That is, the underlying
space   $\C$  is a compact Riemann surface with a finite number of
``holes''. The holes can be ``big'' (e.g. $\mathbb{D}$ can be
thought of as the Riemann sphere with a big hole), or small (e.g.
$\mathbb{C}$ can be thought of as the Riemann sphere with a small
hole).

We begin by stating the main results, postponing the proofs to
\ref{SS:proofs}.

\begin{prop}{\label{P:simpcntd}}
 Every orbifold curve has a universal cover which is a simply connected
 orbifold curve. The simply connected orbifold curves are precisely $\mathbb{H}$,
 $\mathbb{C}$, and $\PP(m,n)$ for $(m,n)=1$.
\end{prop}

\begin{proof}
  Follows from Proposition \ref{P:uniformization}
\end{proof}

\begin{defn}{\label{D:simpcntd}}
 We call an orbifold curve {\bf hyperbolic}, {\bf Euclidean}, or
 {\bf spherical} if its universal cover is $\mathbb{H}$,
 $\mathbb{C}$, or $\PP(m,n)$, respectively.
\end{defn}

 \begin{prop}[{\bf hyperbolic orbifolds}]{\label{P:hyperbolic}}
  All non-algebraic orbifold curves are hyperbolic. An algebraic
  orbifold curve is hyperbolic if and only if
  $$2g-2+\sum\frac{n_i-1}{n_i}$$
  is greater than 0, where   $n_i \leq \infty$ are
     the orders of the orbifold points of ${\mathcal C}$
     (the orbifold order of a hole is $n=\infty$).

  Every hyperbolic orbifold is the quotient of ${\mathbb H}$ by a
  finitely generated discrete subgroup of   $\PSL_2({\mathbb R})$.
  Conversely, the quotient of ${\mathbb H}$ by any such subgroup is a
  hyperbolic orbifold curve.  Two orbifold curves are isomorphic if and
  only if the corresponding subgroups of $\PSL_2({\mathbb R})$ are
  conjugate.
 \end{prop}

\begin{proof}
  See Proposition \ref{P:hEs} and the discussion afterwards.
\end{proof}

 \begin{prop}[{\bf Euclidean orbifolds}]{\label{P:Euclidean}}
   Besides the non-compact ${\mathcal C}_n$ and ${\mathcal C}_{2,2}$, and
   the non-orbifold genus 1 curves, the rest of the Euclidian orbifolds have coarse moduli space
   ${\mathbb P}^1$ and orbifold point orders $(3,3,3)$, $(2,4,4)$,
   $(2,3,6)$ or $(2,2,2,2)$.  The fundamental groups of these are
   $\Lambda_0\rtimes\mu_3$, $\Lambda_0\rtimes\mu_6$,
   $\Lambda_{1728}\rtimes \mu_4$ or $\Lambda\rtimes\mu_2$, which act via the
   obvious embeddings into ${\mathbb C}\rtimes{\mathbb C}^\ast$ on
   ${\mathbb C}$. Here $\Lambda$ denotes a generic
   lattice, and $\Lambda_0=<1,\xi_6>$,
   $\Lambda_{1728}=<1,i>$ are the two symmetric lattices.
 \end{prop}

 \begin{prop}[{\bf spherical orbifolds}]{\label{P:spherical}}
   These all have coarse   moduli space ${\mathbb P}^1$. We have the
   following possibilities:
  \begin{itemize}

   \item[$\mathbf{i.}$] {\em The number of orbifold points is at most 2.}
   In this case $\C\cong\FF(m,n)$, where $n, m\geq 1$ are the orders of the
     orbifold points. We have
     $\pi_1{\mathcal C}\cong{\mathbb Z}_d$, where $d=(m,n)$, and the universal
     cover is ${\mathcal P}(\frac{n}{d},\frac{m}{d})$ on which
     ${\mathbb Z}_d$ acts by rotations.

   \item[$\mathbf{ii.}$] {\em The number of orbifold points is at least 3.}
     In this case we have the following possibilities for the orders
     of the orbifold points:
     $(2,2,n)$, for some $n\geq2$, $(2,3,3)$, $(2,3,4)$ or $(2,3,5)$.
     The corresponding fundamental
     group is the dihedral, tetrahedral, octahedral or icosahedral group,
     respectively,   and the
     universal cover is ${\mathbb P}^1$, acted upon in the usual way.

  \end{itemize}
\end{prop}

\subsection{Proofs}{\label{SS:proofs}}

We begin with computing the fundamental group of an arbitrary
orbifold curve. Let $\C$ be an orbifold curve. Let
$P_1,\cdots,P_k$ be its orbifold points, with
$n_1,\cdots,n_k<\infty$ their  orders, and let $Q_1,\cdots,Q_l$ be
the holes. Assume the genus of the underlying (compactified)
Riemann surface is $g$, and $l+k \geq 1$.

If we remove  the $P_i$ we obtain an open Riemann surface $U$
whose fundamental group is isomorphic to
$\mathbb{Z}^{*(2g+k+l-1)}$. More precisely, let
$a_1,b_1,\cdots,a_g,b_g$ be the generators of $\pi_1\Cm$, $\rho_j$
loops around $P_j$, and $\sigma_h$ loops around $Q_h$. Then,
 {\small
   $$\pi_1U \cong (*\mathbb{Z}a_i)*(*\mathbb{Z}b_i)*(*\mathbb{Z}\rho_j)
          *(*\mathbb{Z}\sigma_h)/<\prod[a_i,b_i]\prod\rho_j\prod\sigma_h=1>.$$}
If we stick the $P_j$ back on, we introduce the new relations
$(\rho_j)^{n_j}=1$, for $1\leq j\leq k$. From van Kampen we obtain
the following

\begin{prop}{\label{P:fundgp}}
   There is an isomorphism
   {\small
   $$\pi_1\C \cong (*\mathbb{Z}a_i)*(*\mathbb{Z}b_i)*(*\mathbb{Z}_{n_j}\rho_j)
          *(*\mathbb{Z}\sigma_i)/<
                            \prod[a_i,b_i]\prod\rho_j\prod\sigma_h=1>.$$ }
\end{prop}

\begin{prop}{\label{P:uniformization}}
  Let $\C$ be an orbifold curve. If $\C$ is a football $\FF(m,n)$
  with $m\neq n$, then the universal cover of
  $\C$ is $\FF(\frac{m}{d},\frac{n}{d})$, where $d=\gcd(m,n)$. In
  particular, $\C$ is not uniformizable. In the remaining cases, $\C$
  is uniformizable. More strongly, we have
  $\C\cong[X/G]$, where $X$ is a Riemann surface, and $G$ is finite group
  acting  on it.
\end{prop}

\begin{proof}
  Let the notation be as in the paragraph preceding the
  proposition.
  We consider several cases  (A-E below). In cases A-D
  we show that the injectivity condition of Theorem \ref{T:uniformizable}
  is satisfied,
  so $\C$ is uniformizable. In fact, we show   more strongly that
  this injectivity condition
  is satisfied for the profinite completion of the $\pi_1\C$.
  So, in cases A-D, $\C$ is in fact the quotient stack of a finite
  group acting on a Riemann surface $C$. Note that,
  if $\C$ is algebraic in any of these cases, then  $C$
  will be algebraic too.

   Case E is the only case in which $\C$ is not uniformizable.
  In this case the universal cover is shown to be $\PP(m,n)$ for
  some relatively prime natural numbers $m$,$n$.

\vspace{0.1in}
 \noindent {\em Case A: $g \geq 1$.} By passing to a quotient,
  we can reduce to the case $k=1$, $l=0$.
  We show that there is a finite quotient
of the group
        $$A=\mathbb{Z}a*\mathbb{Z}b*\Zn c/<aba^{-1}b^{-1}=c>.$$
\noindent into which $\Zn c$ injects. We construct this group as
follows. Consider the action of $\Zn x$ on $\Zn z\oplus\Zn y$
where the action of the generator $x \in \Zn x$ is given by the
matrix
    $$\begin{pmatrix}
         1 & 1 \\
          0 & 1
      \end{pmatrix}$$

The semi-direct product $(\Zn z\oplus\Zn y) \rtimes \Zn x$, the
Heisenberg group over $\Zn$, is made into a quotient of $A$ by
sending $a$,$b$ and $c$ to $x$, $y$ and $z$, respectively. It is
easy to check that this quotient has the desired property.

\vspace{0.1in} \noindent{\em Case B: $g=0$,  $l \geq 1$.} If $k=0$
the statement is trivial. The case $k \geq 1$ is reduced to the
case $k=l=1$ by passing to a quotient, in which case the statement
is again obvious.

\vspace{0.1in} \noindent{\em Case C: $g=0$,  $l=0$ and $k \geq
3$.} Again, by passing to a quotient, we may assume $k=3$. For any
triple of integers $m$, $n$, $p \geq 2$, Fox \cite{Fox} constructs
two permutations $a$, $b$ in a certain finite permutation group
such that $a$ has order $m$, $b$ has order $n$, and $ab$ has order
$p$. This can be used to construct the desired quotient (also see
Remark \ref{R:nonFox}). \vspace{0.1in}

\noindent{\em Case D: $g=0$, $l=0$, and $k=2$;  equal orbifold
orders.} Let $n$ be the orbifold order of the two orbifold points.
Then $\C$ is isomorphic to the quotient of $\mathbb{P}^1$ under
the rotation action of $\Zn$. \vspace{0.1in}

\noindent{\em Case E: $g=0$, $l=0$, and $k=1,2$; unequal orbifold
orders.} When $k=1$, the resulting orbifold, a drop, is  simply
connected (more generally, every $\PP(m,n)$ is simply connected).
The last case is $k=2$. Assume the orbifold points have order
$m\neq n$. Then, the fundamental group is isomorphic to
$\mathbb{Z}_d$, where $d$ is the greatest common divisor of $m$
and $n$. This easily follows from van Kampen. The universal cover
of $\FF(m,n)$ is  $\PP(\frac{m}{d},\frac{n}{d})$; the action of
$\mathbb{Z}_d$ on  $\PP(\frac{m}{d},\frac{n}{d})$ is by
$\frac{2k\pi}{d}$ rotations that fix the two orbifold points.
\end{proof}

\begin{rem}{\label{R:nonFox}}
    There is an alternative method to prove the last case
    ($g=0$, $l=0$, and $k=1,2$)
    of the above proposition. It can be shown that, for any triple
    integers $(m,n,p)$ greater than or equal to 2, and for any
    algebraically closed field $K$ whose characteristic does not divide
    these integers, there are elements $x,y$ in
    $\PSL(2,K)$ such that $x$ has order $m$, $y$ has order $n$, and $xy$ has
    order $p$. The group generated by $x$ and $y$ will be the desired finite
    quotient  of our group.

    To prove the existence of $x$ and $y$, we think of them as automorphisms
    of $\mathbb{P}_K^1$. We take $x$ to be the multiplication
    by a primitive $m$-the root of unity (a rotation fixing $0$ and $\infty$),
    and $y$ to be a rotation of order $n$ fixing $1$ and a variable point
    $v \in \mathbb{P}_K^2$. The condition that the matrix corresponding
    to $xy$ has order $p$ will be a single non-trivial  equation
    in $v$ which can be stated by requiring that the ratio of the
    eigenvalues of the matrix $xy$ is a primitive $p$-root of unity. This
    can always be solved.
 \end{rem}

 Proposition \ref{P:uniformization} immediately implies
 Proposition \ref{P:simpcntd}.
 Proposition  \ref{P:simpcntd} can also be proved
   directly. We outline the proof. Assume $\C$ is simply
   connected. Note that
   the coarse moduli space $\Cm$ of $\C$ must be simply connected
   too
   (otherwise a non-trivial covering space of $\Cm$ pulls back to a non-trivial
   covering space of $\C$ which contradict simple connectivity of $\C$);
   so $\Cm$ is isomorphic to $\mathbb{C}$, $\mathbb{H}$, or $\mathbb{P}^1$.
   Therefore, $\C$
   is one of these spaces, with a bunch of orbifold points on it.
   The fundamental group of such an orbifold is easily computed using
   van Kampen, and one checks that the only possibilities for
   $\C$ are $\mathbb{C}$, $\mathbb{H}$ and orbifold weighted projective
   lines $\PP(m,n)$, for $m$ and $n$ relatively prime.

The following corollaries also follow easily from Proposition
\ref{P:uniformization}.

 \begin{cor}{\label{C:quotient}}
   An orbifold curve is a finite group quotient (of a Riemann surface) if
   and only if its universal cover is
   $\mathbb{H}$, $\mathbb{C}$ or $\mathbb{P}^1$ (i.e., if and
   only if it is a discrete group quotient).
 \end{cor}

 The following is the (orientable case of the)  famous Fenchel
 conjecture which was proved by Fox in 1952
 (see pages 143-144 of \cite{LS}).

 \begin{cor}{\label{C:torsionfree}}
  Let $G$ be a group that is isomorphic to the fundamental
  group of an orbifold (e.g., any group
  of the form discussed in Proposition \ref{P:fundgp}).
  Then, $G$ has a subgroup of finite index that is torsion free.
\end{cor}

\begin{proof}
 If $g\geq 1$, or $l \geq 1$, or $k\geq 3$, then $G$ is isomorphic to the
 fundamental group of a uniformizable Deligne-Mumford curve. So by
 Corollary \ref{C:quotient} $G$ has a torsion free subgroup of
 finite index. (Note that the fundamental group of a Riemann surface is torsion free
 as it acts without fixed points on a 1-dimensional complex domain.)
 If $g=l=0$ and  $k \leq 2$, $G$ is the fundamental group of a football,
 which is a finite cyclic group.
\end{proof}

\begin{prop}{\label{P:hEs}}
  Let ${\mathcal C}$ be an orbifold curve with
  underlying analytic curve $C$.
 \begin{itemize}

   \item[$\mathbf{i.}$] If $C$ is compact of genus $g$,
     then ${\mathcal C}$ is hyperbolic (Euclidean, spherical) if and only if
     $$2g-2+\sum\frac{n_i-1}{n_i}$$
     is greater than 0 (equal to 0, less than 0). Here $n_i < \infty$ are
     the orders of the orbifold points of ${\mathcal C}$.

   \item[$\mathbf{ii.}$] If $C$ is not compact, then
      $\mathcal C$ is hyperbolic, unless
     ${\mathcal C}={\mathcal C}_n$, $1\leq n\leq\infty$, or
     ${\mathcal C}={\mathcal C}_{2,2}$, in which cases it is Euclidean.

  \end{itemize}
\end{prop}

\begin{proof}
  If $\C$ is a football $\FF(m,n)$, the assertion is obvious. So we
  may assume $\C$ is uniformizable (Proposition \ref{P:uniformization}).

   First we show that the only non-compact orbifolds whose
   universal cover is $\mathbb{C}$
are $\C_n$, $1 \leq n \leq \infty$, and $\C_{2,2}$. We can argue
as follows. First of all, the trivial action gives
$\C_1=\mathbb{C}$. Let $G$ be a non-trivial discrete group of
automorphism of $\mathbb{C}$. Then $G$ consists of translations
and rotations. If we only have translation, the possibilities (up
to isomorphism) are $G \cong \mathbb{Z}$ and $G \cong
\mathbb{Z}\oplus\mathbb{Z}$. The latter case is excluded since it
gives a compact quotient. The former gives $\C_{\infty}$. If we
only have rotations, then the only possibilities are finite cyclic
groups of rotations around a fixed point. These give $\C_n$, $2
\leq n < \infty$. Now assume $G$ contains both translations and
rotations. The group of all translations in $G$ is isomorphic to
$\mathbb{Z}$, and without loss of generality we can assume it is
$\mathbb{Z}$ itself (i.e. translation by integers). Conjugation by
any rotation in $G$ takes must take $\mathbb{Z}$ to itself. From
this it easily follows that a rotation in $G$ has order two and it
is centered at a point in $\frac{1}{2}\mathbb{Z}$. There are two
orbifold orbits under this action, namely, $\mathbb{Z}$ and
$\frac{1}{2} + \mathbb{Z}$. The quotient is $\C_{2,2}$.

The rest of non-compact orbifolds are uniformized by $\mathbb{H}$.

Now, assume $\C$ is a compact uniformizable orbifold. Choose a
holomorphic differential form $\omega$ on $\C$. We compute the
degree of $\omega$ as follows: if there were no orbifold points,
then of course $\deg(\omega)=2g-2$. But an orbifold pint of order
$n$ has a contribution of $\frac{n-1}{n}$ (see the end Section of
\ref{SS:SD}). Therefore,
        $$\deg(\omega)=2g-2+\sum\frac{n_i-1}{n_i}.$$
        Using the fact that $\C$ is a finite group quotient  (Corollary
\ref{C:quotient}), it is easily seen that if this degree is
positive (resp., zero, negative), then the universal cover is
$\mathbb{H}$ (resp., $\mathbb{P}^1$, $\mathbb{C}$).
\end{proof}

 It is easy to determine when the expression $2g-2+\sum\frac{n_i-1}{n_i}$
 is negative, zero, or positive. Having done that,  Propositions
 \ref{P:hyperbolic}, \ref{P:Euclidean}, and \ref{P:spherical}
 follow immediately from Proposition \ref{P:hEs}.
 In Proposition \ref{P:hyperbolic}, the fact that the quotient
 of $\mathbb{H}$ under the action of any finitely generated
 discrete subgroup is on orbifold of finite type is not quite trivial.
 It follows from Ahlfors Finiteness Theorem (see \cite{Maskit}).

\vspace{0.2in}
\subsection{Digression on group presentations}

We have shown that the fundamental group of a compact orbifold
whose underlying space is a Riemann surface of genus $g$ has a
presentation of the following form:

   $$G=\langle\begin{array}{lcl}  a_1, b_1,\cdots a_g, b_g,c_1,\cdots,c_k &  | &
         \begin{array}{l}  c_1^{n_1}=c_2^{n_2}=\cdots =c_k^{n_k}=1    \\
           {[a_1,b_1][a_2,b_2]\cdots[a_g,b_g]c_1c_2\cdots c_k=1} \end{array}
                                                           \end{array}\rangle.$$

Conversely, given such a  presentation we can construct an
orbifold (by attaching copies of $\D_n$s to a bouquet of $2g+k$
circles) whose fundamental group is isomorphic to $G$.

Groups defined by such presentations are called $F$-groups and
have been the subject of intense study since the beginning of the
twentieth century. It is unlikely that our methods give any new
results (as perhaps there is nothing left to prove anyway!) but
they do give very clean, short, and conceptual proofs of certain
well known non-trivial facts about $F$-groups (see for instance
Corollary \ref{C:torsionfree}). We will not pursue this further
here and confine ourselves to a few words on triangle
groups.\footnote{Some authors call these {\em von Dyke groups}.}

From the uniformization theorem, we know for what values of
$g$,$k$,$n_1$,$\cdots$,$n_k$ such groups $G$ are finite. For
instance, when $g=0$ and $k=3$, we recover the finiteness theorems
about triangle groups: it is well-known that for the values
               $$(2,2,n), (2,3,3), (2,3,4), (2,3,5)$$
\noindent these groups are finite (the so-called, spherical
triangle groups). These corresponds to orbifolds whose universal
cover is $\mathbb{P}^1$. For
               $$(3,3,3), (2,4,4), (2,3,6)$$
\noindent we obtain the Euclidean triangle groups which correspond
to the orbifolds whose universal cover is $\mathbb{C}$ (these are
of course infinite). For other values  of $(n_1,n_2,n_3)$ we get
the hyperbolic triangle groups. These ones embed as discrete
subgroups of $\PSL(2,\mathbb{R})$. It is well-known that any
hyperbolic triangle group (or in fact, any infinite $F$-group) has
a subgroup of finite index that is torsion free. This corresponds
to the fact that the corresponding orbifold has a finite \'{e}tale
cover that is a Riemann surface.

In general, any finitely generated discrete subgroup $G$ of
$\PSL(2,\mathbb{R})$ is the fundamental group of an orbifold
(Ahlfors Finiteness Theorem), hence it has a presentation of the
above form. If the fundamental domain is a bounded polygon, the
corresponding orbifold will be compact  and algebraic. If we allow
cusps as vertices of the polygon, the corresponding orbifold will
still be algebraic but non-compact. In other cases our orbifold
will only be analytic. In the latter two cases
 $G$ will indeed be isomorphic to a  finite free product of cyclic
groups (because the fundamental group of a non-compact orbifold
curve is so). Conversely, every finite free product of cyclic
groups occurs as the fundamental group of an orbifold (hence,
embeds as a discrete subgroup of $\PSL(2,\mathbb{R})$).

 \section{Interlude: Gerbes over orbifold weighted projective lines}
      {\label{S:GoverF}}

 To prepare for the uniformization of general Deligne-Mumford
 curves (Section \ref{S:UDMG}), in  this section we study gerbes
 over weighted projective lines. We show that the only simply connected
 Deligne-Mumford curves whose
 coarse moduli space is $\mathbb{P}^1$ and have at most two orbifold
 points are weighted projective lines
 (Corollary \ref{C:wpl}).

 Let $\PP(m,n)$ be an orbifold weighted projective line (so $\gcd(m,n)=1$).
 Let $H$ be a finite  group, and let $A=Z(H)$ be its center.
 Our plan is to study Deligne-Mumford curves $\X$ that are
 $H$-gerbes over $\PP(m,n)$, and compute their fundamental
 groups.  What we will do is to make explicit Giraud's
 classification of gerbes over $\PP(m,n)$ in terms
 of their {\em band}, an element in $H^1(\PP(m,n),\Out(H))$,
 and a cohomology class in $H^2(\PP(m,n),A)$.
 Notice that, since  $\PP(m,n)$ is simply
 connected, any $H$-gerbe over it has a trivial band.
 Therefore,
 the set $\Gamma$ of   isomorphism classes of gerbes with
 a trivialized band is in natural bijection with $H^2(\PP(m,n),A)\cong A$.
  Our  explicit construction of this bijection enables us to
  compute the fundamental groups of these gerbes.
 Throughout this section $m,n$ and $H$ are fixed.

 Let us describe how the map $\alpha \:\Gamma \to A$ is defined.
 Let $U, V \subset \PP(m,n)$
 be the complements of the orbifold points of order $m$ and $n$, respectively.
 Fix a base point $*$ on $W=U \cap V$, and an isomorphism
 $\pi_1W \to \mathbb{Z}$. Let $\X$ be an $H$-gerbe over
 $\PP(m,n)$ with a trivialized band, and let $\U$, $\V$ and $\W$ denote its
 restrictions to
 $U$, $V$ and $W$, respectively. Observe that, since $\W$ has trivial band,
 it is indeed a trivial gerbe over $W$. Fix a lift of $*$ to $\X$
 (the bijection $\Gamma \to A$ will be independent of the choice
 of $*$ and its lift to $\X$).  Pick an element $z$ in the center of
 $\pi_1\W$ which maps to $1 \in \pi_1W=\mathbb{Z}$. Let
 $x$ be the image of the $n^{th}$ power of $z$ in $G=\pi_1\U$, and
 let $y$ be the
 image of the $m^{th}$ power of $z^{-1}$ in  $G'=\pi_1\V$.
 These elements both lie in the center of $H$, namely in $A$. Since everything
 takes place inside the abelian group $A$ we will use the additive notation.
 We define $\alpha(\X):=mx+ny \in A$. If instead of $z$ we pick another element
 $z'$ (necessarily of the form $z+a$, for some
 $a \in A \subseteq H \subseteq \pi_1\W$), we would get the pair
 $x+na$ and $y-ma$, which gives rise to the same element $m(x+na)+n(y-ma)=mx+ny$.
 So the map $\alpha \:\Gamma \to A$ is well-defined.

 Next, we construct the
 inverse $\gamma \: A \to \Gamma$. Fix a pair of integers $r$ and $s$  such that
 $rn+sm=1$. Let $a$ be an element in $A$. Let $\U$ be the gerbe
 over $U$ associated to the class of $sa$ in $A/nA$ (see Section \ref{S:Pre}).
 Similarly, let $\V$ be  the gerbe over $V$
 associated to the class of $ra$ in $A/mA$. We want to glue these gerbes
 together to get a gerbe over $\PP(m,n)$. We have the following
 commutative diagrams:
  {\footnotesize
  $$\xymatrix{ 1\ar[r] & H\ar[r]           & G\ar[r]
                                                & \Zn\ar[r]              & 1 \\
               1\ar[r] & H\ar[r]\ar[u]^{=} & H\times\mathbb{Z}\ar[r]\ar[u]
                                               & \mathbb{Z}\ar[r]\ar[u] & 1 }$$
  }

  {\footnotesize
  $$\xymatrix{ 1\ar[r] & H\ar[r]\ar[d]_{=}  & H\times \mathbb{Z}\ar[d]\ar[r]
                                                 & \mathbb{Z}\ar[r]\ar[d] & 1 \\
             1 \ar[r] & H \ar[r]            & G'\ar[r]
                                              & \Zm\ar[r]              & 1 }$$}
 Here $G:=H \times \mathbb{Z}/(sa,-n)=\pi_1\U$.
 Similarly, $G':=H \times \mathbb{Z}/(ra,-m)=\pi_1\V$.
 The $H \times \mathbb{Z}$ that appears in the bottom row of the first
 diagram is identified with $\pi_1(\U|_W)$, and the one on the top row of
 the second diagram is identified with $\pi_1(\V|_W)$.
 Gluing  $\U$ and $\V$ along $W$ (in a way that the trivialization of the
 bands is preserved) is the same as giving an isomorphism from the bottom
 sequence of the first diagram to the top
 sequence of the second diagram (inducing identity on the two ends). We do this
 via the map that is identity on $H$ and multiplication by -1 on $\mathbb{Z}$.
 This produces the desired $H$-gerbe $\gamma(a)$.

 \vspace{0.1in}
 It is straightforward to check that $\alpha$ and $\gamma$ are inverse to each
 other.
\vspace{0.15in}

 It would be illuminating
 to compare the above construction with the Mayer-Vietoris exact sequence
 relative to the open cover $\{U,V\}$ (cohomologies with coefficients in $A$):
{\scriptsize
$$\xymatrix@=7pt{\cdots H^1(U\cup V)\ar[r]\ar@{=}[d] & H^1(U)\oplus H^1(V)\ar[r]\ar@{=}[d]   &
                        H^1(U\cap V)\ar[r]\ar@{=}[d] & H^2(U\cup V)\ar[r]\ar@{=}[d] &
                        H^2(U)\oplus H^2(V)\ar[r]\ar@{=}[d] & H^2(U\cap V)\cdots\ar@{=}[d] \\
                 0\ar[r]  & A[n]\oplus A[m]\ar[r] & A\ar[r]_{mn} &  A\ar[r]_(.28){(s,r)}
                                                                 & A/nA\oplus A/mA\ar[r]   &   0}$$}
 Now, let $a$ be an element in $A$ and $\X$ its corresponding gerbe.
 Our next claim is that there is a natural isomorphism $\pi_1\X\cong H/a$.
 For this we use the explicit construction of the gerbe $\gamma(a)$ as described
 above. We will apply van Kampen  to the covering $\U, \V$ of $\X$.
The van Kampen diagram looks as follows:
 $$\xymatrix@1@R=3pc@C=0pc{ H\times\mathbb{Z}
      \ar[r]^-{1 \mapsto 1}\ar[d]_{\scriptsize{\rotatebox{-90}{\rotatebox{90}{1} $\mapsto$ \rotatebox{90}{-1}}}}
                                                                                 & H\times\mathbb{Z}/(sa,-n) \\
                               \hspace{0.5in}   H \times \mathbb{Z}/(ra,-m)      &                            }$$
 \noindent It is  easy to  check that the corresponding pushout is naturally isomorphic
 to $H/a$.

\vspace{0.1in}
 Summarizing, we have the following

\begin{prop}{\label{P:GoverF1}}
Let $m,n$ be relatively prime positive integers. Let $H$ be a
finite group, and  $A$ its center. Then, the isomorphism classes
of $H$-gerbes $\X$ over $\PP(m,n)$  with a trivialized band are in
natural bijection with $A$. Furthermore,  we have a natural
isomorphism $\pi_1\X \cong H/a$, where $a \in A$ is the element
corresponding to $\X$.
\end{prop}

The element $a$ has an interpretation in terms of the homotopy
fiber exact sequence of the fibration $BH \to \X \to \PP(m,n)$:
    $$\xymatrix@=7pt{\cdots \pi_2BH\ar[r]\ar@{=}[d] & \pi_2\X\ar[r]\ar@{=}[d]  &
                      \pi_2\PP(m,n)\ar[r]\ar@{=}[d] & \pi_1BH\ar[r]\ar@{=}[d]    &
                            \pi_1\X\ar[r]\ar@{=}[d] & \pi_1\PP(m,n)\cdots\ar@{=}[d] \\
                      0\ar[r]  & \mathbb{Z}\ar[r]_d & \mathbb{Z}\ar[r]_{1 \mapsto a}&  H\ar[r]
                                                                      & H/a\ar[r]   &   0}$$
\noindent Here, $d$ is the order of $a$ in $H$. We also see from
the above exact sequence that $\pi_2\X$ is naturally isomorphic to
$\mathbb{Z}$ and maps to $d\mathbb{Z} \subseteq
\mathbb{Z}=\pi_2\PP(m,n)$. The higher homotopy groups $\pi_i(\X)$,
$i\geq 3$, map isomorphically to $\pi_i\PP(m,n)$, which are
themselves isomorphic to $\pi_iS^2$. (To prove the last statement,
use the homotopy fiber exacts sequence for $\mathbb{C}^2-\{0\} \to
\PP(m,n)$).

So far all the gerbes under consideration came with a fixed
trivialization of their band. However, we are only interested in
classification of the underlying stack of these gerbes. To this
end, we have to mod out by the action of $\Out(H)$ on $A$.

\begin{prop}{\label{P:GoverF2}}
Let $m,n$ be relatively prime positive integers. Let $H$ be a
finite group. Then there is a natural bijection
$$\{\text{Algebraic stacks $\X$ that are $H$-gerbes over
$\PP(m,n)$}\}_{/\text{iso}} \rightleftharpoons
{\textcolor{white}.}_{\Out(H)}\backslash A.$$ The fundamental
group of the stack $\X_a$ corresponding to an element $a \in A$ is
(non-canonically) isomorphic to $H/a$.
\end{prop}

 By the above proposition, if we want our gerbe to be simply connected,
 we have to assume $H$ is cyclic and $a$ is a generator for it.
 On the other hand, the set of generators of a cyclic group $H$ is
 acted on transitively by $\Out(H)$. Therefore, for a fixed cyclic
 group $H \cong \Zd$ there is a unique Deligne-Mumford gerbe that is an
 $H$-gerbe over the weighted projective line $\PP(m,n)$. It is nothing but the
 weighted projective line $\PP(md,nd)$ (see Section {\ref{S:Pre}}).

 \begin{cor}{\label{C:wpl}}
 Let $\X$ be a Deligne-Mumford curve whose coarse moduli space is $\mathbb{P}^1$
 and has at most two orbifold points. Then, $\X$ is simply connected
 if and only if is isomorphic to a weighted projective line $\PP(m,n)$.
 \end{cor}

\section{Uniformization  of Deligne-Mumford curves: the general case}
{\label{S:UDMG}}

 In this section we extend the classification of Deligne-Mumford orbifolds
 by their uniformization type (Section \ref{S:UDMO}) to all Deligne-Mumford
 curves. We saw in the previous section that in this case we have a bigger
 collection of simply connected Deligne-Mumford curves. Namely, we have to bring
 into play all $\PP(m,n)$ for $m,n \geq 1$. In fact, it turns out that there
 is no more to add (Proposition \ref{P:simpcntd2}).

 Let $\X$ be an arbitrary Deligne-Mumford curve. Then $\X$ is an $H$-gerbe over
 an orbifold $\C$.
 The orbifold $\C$ is easy to describe: it has the same
 coarse moduli space as $\X$, and at any orbifold point of $\X$ where
 the inertia degree jumps, say by $n$, it has an orbifold point of order $n$.

To classify simply connected Deligne-Mumford curves, we need the
following simple lemma.

\begin{lem}{\label{P:simplyconnected}}
     Let $\X$ and $\Y$ be connected topological stacks, and assume $\X$ is
     an $H$-gerbe over $\Y$ for some discrete group $H$. Assume
     $\X$ is simply connected. Then so is $\Y$.
\end{lem}

\begin{proof}
  This follows immediately from the homotopy fiber sequence for
  the fibration $\X \to \Y$. Alternatively, let $\tilde{\Y}$ be
  the universal cover of $\Y$ and set $\X':=\tilde{\Y}\times_{\Y}\X$.
  Then $\X' \to \X$ is a covering map between connected stacks,
  with $\X$ simply connected, so  it is an isomorphism.
  This implies that $\tilde{\Y} \to \Y$ is also an isomorphism.
\end{proof}

\begin{prop}{\label{P:simpcntd2}}
  Every Deligne-Mumford curve has a universal cover which is a simply connected
Deligne-Mumford curve.  The simply connected Deligne-Mumford
curves are precisely ${\mathbb H}$, ${\mathbb C}$ and ${\mathcal
P}(m,n)$, for arbitrary $n,m\geq1$.
\end{prop}

\begin{proof}
   Assume $\X$ is an $H$-gerbe over an orbifold $\C$.
   By Proposition \ref{P:simplyconnected}, $\C$ is simply
   connected. If $\C$ is $\mathbb{C}$ or $\mathbb{H}$, then
   $\X\cong\C$, as there is no non-trivial gerbe over a
   contractible space. Otherwise, $\C$ is a football, in which
   case the result follows from Corollary \ref{C:wpl}.
\end{proof}

 \begin{defn}
  We call a Deligne-Mumford curve {\bf
  hyperbolic}, {\bf Euclidean} or {\bf spherical}, if its universal
 cover is ${\mathbb H}$, ${\mathbb C}$ or a ${\mathcal P}(m,n)$, respectively.
 \end{defn}

 \begin{prop}{\label{P:type}}
   Let ${\mathcal X}$ be a Deligne-Mumford curve with underlying orbifold
   curve ${\mathcal C}$.  Then ${\mathcal X}$ is hyperbolic (Euclidean,
   spherical) if and only if ${\mathcal C}$ is so.
 \end{prop}

 \begin{proof}
   The case where $\C$ is hyperbolic or Euclidean follows from Proposition
   \ref{P:characterize}, as in this case $\pi_n\C$ vanishes for $n \geq 2$.
   When $\C$ is spherical, the composite map
   $\tilde{\X} \to \X \to \C$
   induces a surjective map
       $$\tilde{\X}_{mod} \to \Cm\cong\mathbb{P}^1$$
   on the coarse moduli spaces.
   So, $\tilde{\X}$ can not be any of  $\mathbb{C}$ or $\mathbb{H}$.
   Therefore, $\X$ is spherical.
 \end{proof}

 \begin{prop}[{\bf hyperbolic and Euclidean curves}]{\label{P:hE}}
  Assume $\X$ is a hyperbolic or Euclidean Deligne-Mumford analytic curve.
  Let $\C$ be the orbifold over which $\X$ is a gerbe, and let $H$ be the
  generic  inertia group of $\X$.  We have the short exact sequence
  $$ 1 \to H \to \pi_1\X \to \pi_1\C \to 1.$$
  Furthermore, the orbifold $\C$ together with the above exact sequence
  uniquely determine $\X$ up to isomorphism. More precisely,
  we have $\X\cong[\tilde{\C}/\pi_1\X]$, where $\pi_1\X$ acts on
  $\tilde{\C}$ via $\pi_1\C$.
 \end{prop}

 \begin{proof}
   Follows from Proposition \ref{P:characterize}.
 \end{proof}

 \begin{prop}[{\bf spherical curves}]{\label{P:spherical2}}
  Assume $\X$ is a spherical Deligne-Mumford analytic curve.
  Let $\C$ be the orbifold
  over which $\X$ is a gerbe, and let $H$ be the generic
  inertia group of $\X$.
   \begin{itemize}

     \item[$\mathbf{i.}$] $\ker(H \to \pi_1\X)$ is cyclic, say of order $d$,
        and we have an exact sequence
           $$ 1 \to \Zd \to H \to \pi_1\X \to \pi_1\C \to 1.$$

     \item[$\mathbf{ii.}$] If $\C$ has orbifold type $(m,n)$, then the universal
        cover of $\X$ is isomorphic to $\PP(d\frac{m}{a},k\frac{n}{a})$, where
        $a=\gcd(m,n)$.

     \item[$\mathbf{iii.}$]  If $\C$ has orbifold type $(2,2,n)$, $(2,3,3)$,
       $(2,3,4)$ or $(2,3,5)$, then the universal cover of $\X$ is isomorphic
       to $\PP(d,d)$.

\end{itemize}
 \end{prop}

\begin{proof}
  The exact sequence is simply the fiber homotopy exact sequence
  of the fibration $\X \to \C$ (note that $\pi_2\C\cong\mathbb{Z}$).
  This  proves Part ($\mathbf{i}$).

 \vspace{0.1in}
  \noindent{\em Part} ($\mathbf{ii}$).
  Let $K$ be the kernel of $H \to \pi_1\X$ (so $K \cong \Zd$), and
  let $d=\gcd(m,n)$.
  \noindent Let $\tilde{\X}$ be the universal cover of $\X$.
  We have to show that $\tilde\X \cong \PP(d\frac{m}{a},d\frac{n}{a})$. We know
  that $\tilde{\X}$ is an $H'$-gerbe over a simply connected orbifold
  $\tilde{\C}$ which maps via a covering map to
  $\C$ (hence the notation $\tilde{\C}$!). So we have
  $\tilde{\C} \cong \PP(\frac{m}{a},\frac{n}{a})$. On the other hand,
  we know by Lemma \ref{L:hidden} that the kernel of $H' \to \pi_1\tilde{\X}=1$ is
  the same as the kernel of $H \to \pi_1\X$; hence, $H'\cong \Zd$.
  Therefore, $\tilde\X \cong \PP(d\frac{m}{a},d\frac{n}{a})$, by
  Corollary \ref{C:wpl}.

 \vspace{0.1in}
  \noindent{\em Part} ($\mathbf{iii}$).
  Using an argument similar to the previous part, it follows that the universal
  cover of $\X$ is the (unique) simply connected $\Zd$-gerbe over $\mathbb{P}^1$
  This is isomorphic to $\PP(d,d)$.
 \end{proof}

 \begin{cor}{\label{C:quotient2}}
   Let $\X$ be a uniformizable Deligne-Mumford curve. Then there exists
   a Riemann surface $X$ and a finite group $G$ such that $\X\cong[X/G]$.
 \end{cor}

 \begin{proof}
   When $\X$ is spherical the result is obvious. So we assume $\X$ is not
   spherical. It is enough to show that there exists a Riemann surface $X$
   with a finite covering map $X \to \X$. Assume $\X$ is an $H$-gerbe over
   an orbifold $\C$. It is enough to show that there exists a finite
   covering map $X \to \C$, with $X$ Riemann surface, such that the
   induced $H$ gerbe on $X$ is trivial.

   By Proposition \ref{P:type}, $\C$ is uniformizable. So, by Corollary
   \ref{C:quotient}, there is a finite covering map $C \to \C$.
   So, by  pulling back our $H$-gerbe to $C$, we are reduced to the case
   wehere $\C=C$ is a Riemann surface. Recall now that an $H$-gerbe
   over $C$ is characterized by its band, which is an element
   $\beta\in H^1(C,\Out(H))\cong\Hom(\pi_1C,\Out(H))$, and an element in
   $\alpha\in H^2(C,Z(H))$. Observe that $\beta$ is the
   homomorphism $\beta \: \pi_1C \to \Out(H)$ induced from the exact sequence
   of Proposition \ref{P:hE}. We can easily kill $\beta$ by passing to
   a finite cover of $C$. (For instance, let $K$ be the kernel
   of the homomorphism $\pi_1C \to \Out(H)$, and replace $C$ by the covering
   space $X \to C$ corresponding to the subgroup $K \subseteq \pi_1C$.)

   To show that we can kill $\alpha$, we show that, in general, for any finite
   abelian group $A$, any element $\alpha \in H^2(C,A)$ can be killed after
   passing to a finite cover $f \:X \to C$. If $C$ is not compact, $H^2(C,A)$
   is already trivial and there is nothing prove. So assume $C$ is compact.
   It is an easy exercise
   that $H^2(C,A)\cong A$ and the map $f^* \: H^2(C,A) \to H^2(X,A)$ is
   simply multiplication by $d$, where $d$ is the degree of $f$. Since
   $\alpha \in H^2(C,A)$ is torsion, say $n\alpha=0$, it is enough
   find such an $f$ of degree $n$. Equivalently, we have to show
   that the fundamental group
   $\pi_1C\cong(*\mathbb{Z}a_i)*(*\mathbb{Z}b_i)/<\prod[a_i,b_i]=1>$
   has a subgroup of index  $n$. Note that genus of $C$ is at
   least $1$ (because $\X$, and so $\C$ as well, are not spherical).
   It is easy to see that the subgroup generated by
   $\{a_1,\cdots,a_g,b_1^n,b_2,\cdots,b_g\}$ is such a subgroup.
   This completes the proof.
 \end{proof}

 As we see in Proposition \ref{P:spherical2},   spherical
 Deligne-Mumford curves are
 more difficult to classify than their  Euclidean and hyperbolic
 counterparts: in the hyperbolic and Euclidean cases, the exact sequence of
 Proposition \ref{P:hE} completely classifies the curve in
 question; but in
 the spherical case the exact sequence of Proposition \ref{P:spherical2} is not
 enough to classify the corresponding   Deligne-Mumford curve.
 For instance, part of the information that is missing in the
 exact sequence of  Proposition \ref{P:spherical2} is that
 the homomorphism $H \to \pi_1\X$ carries a natural crossed-module
 structure; this is a general fact about fundamental
 groups of gerbes whose proof will not be given here (but in
 this special case it follows from the results of Section \ref{S:spherical}).
 Unfortunately, even
 the knowledge
 of this crossed-module structure is not enough to determine our Deligne-Mumford
 curve; but it is a step in the right direction. We address this issue
 in the subsequent sections by bringing some 2-group theory into the
 picture.

\subsection{Graphs of groups as skeleta of open Deligne-Mumford
curves}{\label{SS:graphs}}

There is an alternative way of computing the homotopy groups of an
{\em open} Deligne-Mumford curve, which is in a way simpler and
more illuminating. The idea is that any open Deligne-Mumford curve
has the homotopy type of a {\em graph of groups} in the sense of
Serre \cite{Serre}. We do not intend to give an overview of the
theory of graphs of groups here, so we just confine ourselves to
describing the procedure by which we associate a graph of groups
to a Deligne-Mumford curve, and we leave it to the interested
reader to verify that this approach will yield the same results
(e.g., the same formula for the fundamental group etc.).

 The idea is that one can retract an open Deligne-Mumford curve to a
graph of groups. Assume $\C$ is a smooth Deligne-Mumford curve
that is generically an $H$-gerbe, for some finite group $H$. Let
$g$ be the genus of the underlying Riemann surface. Let
$P_1,P_2,\cdots,P_k$, $k \geq 0$, be the orbifold points of $\C$,
with their respective inertia groups $G_1,G_2,\cdots,G_k$,  and
let $Q_1,Q_2,\cdots,Q_l$, $l \geq 1$, be the punctured points.

Let us for one moment forget about the stacky structure of $\C$
and look at the underlying surface. This surface can be deformed
into a bouquet of $2g+l-1$ circles. If we take into account the
stacky structure we should be careful how we treat the orbifold
points $P_1$,$P_2$,$\cdots$,$P_k$. What we do is we choose a base
point $P_0$, different from the $P_i$, and join $P_0$ to each
orbifold point by a path (Figure \ref{fig1}).
 \begin{figure}[!ht]
  \begin{center}
   \input{fig1.pstex_t}
    \caption{}
    \oldlabel{fig1}
  \end{center}
 \end{figure}
\noindent Now we perform the retraction of $\C$ in a way that all
these paths are left intact during the retraction. The result is a
graph of groups whose underlying graph consists of $2g+l-1$
circles together with $k$ segments, all joined at a single point
(Figure \ref{fig2}).
 \begin{figure}[!ht]
  \begin{center}
   \input{fig2.pstex_t}
    \caption{}
    \oldlabel{fig2}
  \end{center}
 \end{figure}
The group associated to each edge is $H$, the group associated to
$P_0$ is also $H$. The group associated to $P_i$, $i \geq 1$, is
$G_i$. For every edge $P_0P_i$ The homomorphism $H \to G_i$ is
just the natural inclusion (obtained from the path $P_0P_i$). For
every circle in the bouquet, the homomorphism $H \to H$ is the one
obtained from the monodromy action  (see footnote in Section
\ref{SS:SD}). This action determines a semi-direct product
              $$H\rtimes (\underset{ 2g+l-1}{*}\mathbb{Z})$$
\noindent which is the fundamental group of the bouquet. The
fundamental group of the whole graph of groups is
         $$(H\rtimes (\underset{ 2g+l-1}{*}\mathbb{Z}))*_H(*_HG_i).$$

In the case of a compact Deligne-Mumford curves one has to use
2-dimensional complexes of groups. The complexes of groups that
arise in this way are of particularly simple form, because they
can be obtained by attaching  a simplex (with generic group $H$)
to the complex of groups of an open Deligne-Mumford curves (which
is homotopy equivalent to a graph of groups). This way  we
recover, for instance, the formula for the fundamental group of an
arbitrary Deligne-Mumford curve (using van Kampen).

 \section{Weighted projective general linear 2-groups}{\label{S:wtd}}

 In this section
 we  invoke the full 2-categorical structure of the category
 of stacks, so our convention in Section \ref{S:NT} that all 2-isomorphic
 morphism are declared equal is no longer in effect.

 As hinted in the previous section (see just after Proposition
 \ref{P:spherical2}), crossed-modules
 (equivalently, {\em 2-groups})  appear naturally in the
 study of fundamental groups of Deligne-Mumford stacks. Even
 in the hyperbolic and Euclidean cases (Proposition \ref{P:hE}) the 2-groups
 are implicitly present: $H \hookrightarrow \pi_1\X$ should be thought of as a
 crossed-module.
 The reason 2-groups occur is that the automorphisms $\Aut\X$ of a stack $\X$
 (viewed as an object in the 2-category of stacks) naturally form a
 weak\footnote{In fact, we still have the associativity on the nose.
 But inverses only exists in the weak sense.} 2-group.
 A quick review of the basic facts of 2-group theory
 is given in Appendix.

 Our main task in this section, and the next, is to study the 2-group
 of automorphisms of weighted projective stacks. Let $n_1,n_2,\cdots,n_r$
 be a sequence of positive integers, and
 consider the weight $(n_1,n_2,\cdots,n_r)$ action of $\mathbb{C}^*$ on
 $\mathbb{C}^r-\{0\}$ (so $t \in \mathbb{C}^*$ acts as multiplication by
 $(t^{n_1},t^{n_2},\cdots,t^{n_r})$). The stacky quotient of this action
 is the weighted projective stack $\PP(n_1,n_2,\cdots,n_r)$. We define
 the {\bf weighted projective general linear 2-group}
 $\PGL(n_1,n_2,\cdots,n_r)$
 to be the 2-group associated to crossed-module
 $[\varphi \: \mathbb{C}^* \to G_{n_1,n_2,\cdots,n_r}]$, where
 $G_{n_1,n_2,\cdots,n_r}$ is the group of all $\mathbb{C}^*$-equivariant
 (for the above weighted
 action) algebraic automorphisms
 $f \: \mathbb{C}^r-\{0\} \to \mathbb{C}^r-\{0\}$.
 The homomorphism $\varphi \: \mathbb{C}^* \to G_{n_1,n_2,\cdots,n_r}$
 is the one induced from the
 $\mathbb{C}^*$-action. We take the action of $G_{n_1,n_2,\cdots,n_r}$
 on $\mathbb{C}^*$ to be
 trivial.

 Throughout this section we fix the weight sequence $n_1,n_2,\cdots,n_r$
 and denote $G_{n_1,n_2,\cdots,n_r}$ simply by $G$.

  \begin{thm}{\label{T:2-aut}} Let $\Aut^{an}\PP(n_1,n_2,\cdots,n_r)$ and
  $\Aut\PP(n_1,n_2,\cdots,n_r)$  be the 2-groups
  of analytic and algebraic automorphisms of $\PP(n_1,n_2,\cdots,n_r)$,
  respectively. Then,
  \begin{itemize}

   \item[$\mathbf{i.}$] The analyfication functor  induces an equivalence
      $$\Aut\PP(n_1,n_2,\cdots,n_r) \to \Aut^{an}\PP(n_1,n_2,\cdots,n_r)$$
      of 2-groups.
   \item[$\mathbf{ii.}$] The natural map
     $$ \PGL(n_1,n_2,\cdots,n_r) \to \Aut\PP(n_1,n_2,\cdots,n_r)$$
     is an equivalence of 2-groups.
  \end{itemize}
 \end{thm}

 First we prove a

  \begin{lem}
    Let $H$ be an abelian group scheme acting on a connected scheme $X$ over
    $\mathbb{C}$, and let $\X=[X/H]$ be the quotient stack.
    Let  $G$ be the group of $H$-equivariant automorphism of $X$.
    Let $\Aut\X$ be the (weak) 2-group of automorphism of $\X$.
     \begin{itemize}

    \item[$\mathbf{i.}$] The natural homomorphism
       $\varphi \: H(\mathbb{C}) \to G(\mathbb{C})$ can be turned into
       a crossed-module by taking the trivial action of $G$ on $H$.

    \item[$\mathbf{ii.}$] Let $\mathfrak{G}$ denote the
      2-group associated to the crossed-module $[\varphi \: H \to G]$.
      Then there is a natural  map of 2-groups
      $\mathfrak{G} \to \Aut\X$. Furthermore, this morphism induces an
      isomorphism on $\pi_2$.

    \item[$\mathbf{iii.}$] Assume that $\X$ is a proper Deligne-Mumford stack
    and $H$ is   affine.
    Then the induces map on $\pi_1$ is injective.

 \end{itemize}

 \end{lem}

  \begin{proof}
      Part ($\mathbf{i}$) is straightforward, because $\varphi$ maps $H$
      into the center of $G$.

  \vspace{0.1in}

  To prove part ($\mathbf{ii}$), we recall the explicit description
      of the quotient stack $[X/H]$:

      {\small
      $$\Ob[X/H](S)=\left\{
      \begin{array}{rcl} (T,\alpha) &\vert& T \  \text{an $H$-torsor over $S$} \\
                             & & \alpha \: T \to X  \ \text{an $H$-map}
                                              \end{array}\right\}$$
      $$\Mor[X/H](S)((T,\alpha),(T',\alpha'))=\{f \: T \to T'
        \ \text{an $H$-torsor map s.t.} \ \alpha'\circ f=\alpha\}$$}

  Any element of $g \in G(\mathbb{C})$ induces
  an automorphism of $\X$ (keep the same
  torsor $T$ and compose $\alpha$ with the action of $g$ on $X$).
    Also, for any element $h \in H(\mathbb{C})$, there is a natural
  2-isomorphism from the identity automorphism of $\X$ to the automorphism
  induced by  $\varphi(h) \in G(\mathbb{C})$ (which is by
  definition the same as the action
  of $h$). It is given by the multiplication
  action of $h^{-1}$ on
  the torsor $T$ (remember $H$ is abelian)
  which makes the following triangle commute
       $$\xymatrix{ T \ar[r]^{\alpha} \ar[d]_{h^{-1}} &  X  \\
                    T \ar[ru]_{h\circ\alpha}    &           }$$

  Interpreted in the language of 2-groups, this gives a
   morphism of 2-groups   $\mathfrak{G} \to \Aut\X$.

  A general fact about stacks is that, the group of 2-isomorphisms
  from the identity automorphism of $\X$ to itself is naturally isomorphic
  to the group of global sections of the inertia stack of $\X$. In the
  case $\X=[X/H]$, this is naturally isomorphic to the  group of elements
  of $H(\mathbb{C})$ that act trivially on $X$.
  Note that this group is naturally
  isomorphic to $\pi_2\mathfrak{G}$. Therefore, the map
  $\mathfrak{G} \to \Aut\X$ induces an isomorphism on $\pi_2$. This completes
  the proof of ($\mathbf{ii}$).

  To prove part ($\mathbf{iii}$) we need to show that, if the action of
  $g \in G(\mathbb{C})$ on $\X$ is 2-isomorphic to identity,
  then $g$ is in the image of
  $\varphi$. Let us fix such a 2-isomorphism. By looking at
  the groupoid $[X/H](X)$
  and following around the effect of the action of $g$  and   our
  2-isomorphism on the point $X \to [X/H]$, we
  end up with an $H$-torsor map $F \: H\times X \to H\times X$ making
  the following $H$-equivariant triangle commute:
  $$\xymatrix{ H\times X \ar[r]^(0.62){\mu} \ar[d]_{F} &  X  \\
                    H\times X \ar[ru]_(0.55){g\circ\mu}    &           }$$
  Here, $H\times X$ is the trivial $H$-torsor on $X$
  and $\mu$ corresponds to the action of $H$ on $X$.

   Composing $F$ with the identity section $X \to H\times X$, we obtain
   a map $f \: X \to H$. It follows from the $H$-equivariance of the
   above diagram that, for any point $x$ in $X$, the effect of $g$ on $x$
   is the same as the effect of $f(x)$ on $x$. In other words, $f(x)g^{-1}$
   leaves $x$, or in fact every point $h(x)$ in the $H$-orbit of $x$, invariant.
   The same is true for $f(h(x))g^{-1}$, for any $h \in H(\mathbb{C})$.
   This implies that,
   for any point $x$ of $X$, and any $h \in H(\mathbb{C})$, the element
   $r(x,h):=f(h(x))f(x)^{-1}$ leaves $x$ (in fact the entire orbit of $x$)
   fixed; that is, $r(x,h)$ belongs to the stabilizer group
   $S_x \subseteq H$ of $x$.
   Therefore,  the map  $\rho \:H\times X \to H\times X$,
   $\rho(h,x):=(r(h,x),x)$  factors through the stabilizer
   group scheme $S \to X$. Since $S \to X$  has discrete fibers,
   the restriction of $\rho$ to $H_0 \times X$
   factors through the identity section ($H_0$ is the connected
   component of the identity).
   Hence, for every $h \in H_0(\mathbb{C})$ and $x \in X$, $r(h,x)$ is the
   identity.
   This implies that $f \: X \to H$ is $H_0$-equivariant
   (for the trivial action of
   $H_0$ on $H$). So, we obtain
   an induced map $[X/H_0] \to H$. But we know that $[X/H_0]$, being
   finite over $[X/H]$, is proper and $H$
   is affine, so this map must be constant; hence $f \: X \to H$ is also
   constant. Let $h \in H(\mathbb{C})$ be the value of $f$. Then, the action
   of $h$ on $X$ coincides with the action of $g$. So  $\varphi(h)=g$, which
   is what we wanted to prove.
  \end{proof}

  \begin{rem}\end{rem}

  \vspace{0.1in}
   \begin{itemize}

    \item[$\mathbf{1.}$] In fact $\Aut(\X)$ is a 2-group scheme
     and the above proof can be easily jazzed up to the level of
     2-group schemes.

    \item[$\mathbf{2.}$]  The assumption on schemes being
    over $\mathbb{C}$ is not necessary; we can work over any base scheme.
    Part ($\mathbf{i}$) and  ($\mathbf{ii}$) of the lemma
    are quite formal and are valid in the analytic and topological settings
    as well. Part ($\mathbf{iii}$) of the lemma is also valid in the analytic
    setting (but we have to require that $[X/H^0]$ is proper, where
    $H^0$ is the connected component of the identity).

    \item[$\mathbf{3.}$] Part ($\mathbf{iii}$) of the lemma can be interpreted
    as saying that the morphism $\mathfrak{G} \to \Aut\X$ is ``injective''.
    That is, its kernel is equivalent to the trivial 2-group.
 \end{itemize}

 \vspace{0.1in}

  \begin{proof}[Proof of Theorem \ref{T:2-aut}.]
     Part ($\mathbf{i}$) is a GAGA type statement and
     follows from Theorem 1.1 of \cite{Lurie}.
     To prove part ($\mathbf{ii}$)
     we may assume $r \geq 2$. By the previous lemma, the map
     $$PGL(n_1,n_2,\cdots,n_r) \to \Aut\PP(n_1,n_2,\cdots,n_r)$$
     \noindent induces an
     isomorphism on $\pi_2$ and an injection on $\pi_1$. So we only
     need to prove the surjectivity on $\pi_1$. First, we show that for
     any automorphism
     $f \: \PP(n_1,n_2,\cdots,n_r) \to  \PP(n_1,n_2,\cdots,n_r)$  there is
     a lift
     $F \: \mathbb{C}^r-\{0\} \to \mathbb{C}^r-\{0\}$ making
     the following diagram 2-cartesian

       $$\xymatrix{\mathbb{C}^r-\{0\} \ar[r]^F \ar[d]
                           & \mathbb{C}^r-\{0\} \ar[d]  \\
               \PP(n_1,n_2,\cdots,n_r) \ar[r]_f  & \PP(n_1,n_2,\cdots,n_r)  }$$

     Recall that the map  $\mathbb{C}^r-\{0\} \to \Pic\PP(n_1,n_2,\cdots,n_r)$
     is the total space of the Hopf bundle
     $L \in \Pic\PP(n_1,n_2,\cdots,n_r)$ with the identity section removed.
     We know that $\Pic\PP(n_1,n_2,\cdots,n_r)\cong\mathbb{Z}$, with the
     Hopf     bundle      corresponding to $1 \in \mathbb{Z}$. Any automorphism
       induces a trivial action on $\Pic\PP(n_1,n_2,\cdots,n_r)$ (i.e. it
       preserves the ``degree'').
     Hence, $f^*L\cong L$.  A choice of an isomorphism
     $f^*L \risom L$ gives rise to the desired
     2-cartesian diagram.

     If we show that $F$ is $\mathbb{C}^*$-equivariant we are done.
     Since the above diagram is 2-cartesian,
     there is  a groupoid map from
     $\mathbb{C}^*\times \mathbb{C}^r-\{0\}
     \sst{} \mathbb{C}^r-\{0\}$ to itself whose
     corresponding map on the `objects-space' is
     $F \: \mathbb{C}^r-\{0\} \to \mathbb{C}^r-\{0\}$.
     Let us denote the map on the `arrow space' by
     $\Phi \:\mathbb{C}^*\times \mathbb{C}^r-\{0\}
     \to \mathbb{C}^*\times \mathbb{C}^r-\{0\}$.
     Since $\Phi$ commutes with the
     source map (projection) $\mathbb{C}^*\times \mathbb{C}^r-\{0\}
     \to \mathbb{C}^r-\{0\}$,
     for any $x \in \mathbb{C}^r-\{0\}$
     we obtain a map $\Phi_x \: \mathbb{C}^* \to \mathbb{C}^*$
     (the first $\mathbb{C}^*$ is the fiber over $x$
     and the second is the fiber over $F(x)$). A priori, the $\Phi_x$ are not
     group homomorphisms, but since $\Phi$ sends the identity to identity, so
     do $\Phi_x$. It is easy to show that the only algebraic (or analytic)
     isomorphisms  $\mathbb{C}^*  \to \mathbb{C}^*$ that preserve $1$ are
     the identity map $x \mapsto x$ and the
     inversion map $x \mapsto \frac{1}{x}$.
     So each $\Phi_x \: \mathbb{C}^* \to \mathbb{C}^*$ is either identity
     or inversion. By continuity, either all $\Phi_x$ are identity, or
     all inversions. In the first case, it follows that
     $F$ is $\mathbb{C}^*$-equivariant, which is what we want.
     In the second case, it follows that $F$ is twisted
     $\mathbb{C}^*$-equivariant; that
     is, for every $a \in \mathbb{C}^*$, and every
     $x \in \mathbb{C}^r-\{0\}$, we have
     $F(a\cdot x)=\frac{1}{a}\cdot x$.
      We show that this is impossible. First observe that, since $r \geq 2$,
      the map $F\: \mathbb{C}^r-\{0\}  \to \mathbb{C}^r-\{0\} $
      extends to $\mathbb{C}^r$. Consider the point
      $(1,0,\cdots,0) \in \mathbb{C}^r-\{0\}$. If $F$ is twisted
      $\mathbb{C}^*$-equivariant, we  have
      $F(a^{n_1},0,\cdots,0)=a^{-n_1}F(1,0,\cdots,0)$, for every
      $a \in \mathbb{C}^*$. Letting $a$ approach zero, we see that the
      left
      hand side approaches to $F(0,0,\cdots,0)$, which is a well-defined point,
      but the right hand side
      goes to infinity. This  contradiction shows that $F$ can not be
      twisted $\mathbb{C}^*$-equivariant, so it has to be
      $\mathbb{C}^*$-equivariant. The proof is now complete.
 \end{proof}

 \begin{prop}{\label{P:automorphism}}
   Let $n_1,n_2,\cdots,n_r$ be a sequence of positive integers, and let
   $d$ be their greatest common divisor. Then, the weighed projective
   stack $\PP(n_1,n_2,\cdots,n_r)$ is a $\Zd$-gerbe over
   $\PP(\frac{n_1}{d},\frac{n_2}{d},\cdots,\frac{n_r}{d})$. The induced
   map   $$ \PGL(n_1,n_2,\cdots,n_r) \to
   \PGL(\frac{n_1}{d},\frac{n_2}{d},\cdots,\frac{n_r}{d})$$
   of automorphism 2-groups (obtained from the functoriality of the
   ``underlying orbifold'' construction,
   see Proposition \ref{P:gerbe})
   induces an isomorphism on $\pi_1$. That is,
     $$ \pi_1\PGL(n_1,n_2,\cdots,n_r) \risom
   \pi_1\PGL(\frac{n_1}{d},\frac{n_2}{d},\cdots,\frac{n_r}{d})$$
   Note that $\PGL(\frac{n_1}{d},\frac{n_2}{d},\cdots,\frac{n_r}{d})$
   is equivalent to a group, that is,
     $$\PGL(\frac{n_1}{d},\frac{n_2}{d},\cdots,\frac{n_r}{d})
      \risom \pi_1\PGL(\frac{n_1}{d},\frac{n_2}{d},\cdots,\frac{n_r}{d}).$$

 \end{prop}

  \begin{proof}
     Observe that
     $G_{n_1,n_2,\cdots,n_r}=G_{\frac{n_1}{d},\frac{n_2}{d},\cdots,\frac{n_r}{d}}$;
    we denote both by $G$.
    On the level of crossed-modules, our  2-group morphism  look as follows
       $$\xymatrix{  \mathbb{C}^* \ar[r] \ar[d]_{-^d} &  \ar[d]^{id} G  \\
                   \mathbb{C}^* \ar[r]        &     G      }$$
 \noindent (see just before Theorem \ref{T:2-aut} for notation). The result
 is obvious.
 \end{proof}

 \subsection{Explicit description of $\PGL(m,n)$}{\label{SS:explicit}}

 The crossed-module $\PGL(m,n)$ can be described quite explicitly.
 To do so, we first have to determine $G_{m,n}$. We do this by
 considering
 three cases. As before, we denote $G_{m,n}$ simply by $G$.

\vspace{0.1in} \noindent  \framebox[3.7cm][l]{\em The case $m < n,
m\nmid n$.}

\vspace{2mm}

    It is easy to see that an element in $G$ is a map of the form
       $$(x,y) \mapsto (\lambda_1x, \lambda_2y),$$
    for some fix $\lambda_1,\lambda_2 \in \mathbb{C}^*$.
    So we have
        $G \cong \mathbb{C}^*\times \mathbb{C}^*$.
  The map $\varphi \: \mathbb{C}^* \to G$ is given by
    $$\lambda \mapsto (\lambda^m,\lambda^n).$$
 The map
    $$\begin{array}{rcl}
         G \cong \mathbb{C}^*\times \mathbb{C}^* & \to & \mathbb{C}^* \\
               (\lambda_1,\lambda_2)& \mapsto &
                              \lambda_1^{\frac{n}{d}}\lambda_2^{-\frac{m}{d}}
    \end{array}$$
  induces an isomorphism  $\pi_1\PGL(m,n) \cong \mathbb{C}^*$.
  We can produce a section
             $$\sigma \: \pi_1\PGL(m,n) \to G$$
  for this map by choosing integers
  $r$,$s$ such that $sm+rn=d$ and setting
  $\sigma(\lambda)=(\lambda^{-r},\lambda^{s})$.

  \vspace{0.1in}
  \noindent{\em Convention.} Throughout the paper,  given
  $m$, $n$, we fix a pair $r$,$s$ with $sm+rn=d$. When $m$ divides $n$ we take
  $s=1$ and $r=0$.

\vspace{3mm}

\noindent  \framebox[3.7cm][l]{\em The case $m < n, m\mid n$.}

\vspace{2mm}

    Every element in $G$ is a map of the form
       $$(x,y) \mapsto (\lambda_1x, \lambda_2y+ax^{\frac{n}{m}}),$$
    for some $\lambda_1,\lambda_2 \in \mathbb{C}^*$ and $a \in \mathbb{C}$.
   So we have
     $$G \cong (\mathbb{C}^*\times\mathbb{C}^*)\ltimes\mathbb{C},$$
    where the action of an element $(\lambda_1,\lambda_2)
    \in \mathbb{C}^*\times\mathbb{C}^*$
    on an element $a \in \mathbb{C}$ is given by
       $$(\lambda_1,\lambda_2)\cdot a=\lambda_1^{-\frac{n}{m}}\lambda_2a.$$
    (Note the similarity with the group of $2\times 2$ lower-triangular
    matrices.)

 The map $\varphi \: \mathbb{C}^* \to G$ is given by
    $$\lambda \mapsto (\lambda^m,\lambda^n,0).$$
 The map
    $$\begin{array}{rcl}
         G \cong (\mathbb{C}^*\times \mathbb{C}^*)\ltimes\mathbb{C} &
                                  \to & \mathbb{C}^*\ltimes\mathbb{C} \\
               (\lambda_1,\lambda_2,a)& \mapsto &
                                     (\lambda_1^{-\frac{n}{m}}\lambda_2,a)
    \end{array}$$
  induces an isomorphism $\pi_1\PGL(m,n) \cong \mathbb{C}^*\ltimes\mathbb{C}$.
  In the latter semi-direct product the action of $\mathbb{C}^*$ is simply the
  multiplication action.

  There is  a section $\sigma \: \pi_1\PGL(m,n) \to G$ to the above map which is
  given by   $\sigma(\lambda,a)=(1,\lambda,a)$.

\vspace{3mm}

\noindent \framebox[2.7cm]{\em The case $m=n$.}

\vspace{2mm}

We switch the notation from $m=n$ to $d$. In this case we
obviously have $G\cong \GL_2$. The map
    $\varphi \: \mathbb{C}^* \to G$ is given by
      $$\lambda \mapsto \left( \begin{array}{cc} \lambda^d  & 0 \\
                  0  & \lambda^d \end{array}\right) $$

\vspace{2mm} We summarize the above discussion in the following
propositions.

\begin{prop}{\label{P:split}}
   Let $m$, $n$ be distinct positive integers, and let $d=\gcd(m,n)$.
   Then $\PGL(m,n)$ is a split 2-group. In particular, it is classified by its
   homotopy groups:
       $$\begin{array}{rcl}
          \pi_1\PGL(m,n)  &\cong&
                    \left\{  \begin{array}{l}
                        \mathbb{C}^*,  \ \text{if}\ m < n, m\nmid n \\
            \mathbb{C}^*\ltimes\mathbb{C},  \ \text{if}\ m < n, m\mid n
                                                        \end{array}  \right.  \\
           \pi_2\PGL(m,n)  &\cong&    \mathbb{Z}_d.
        \end{array}$$
    (In the case $m\mid n$, the action of $\mathbb{C}^*$ on $\mathbb{C}$
    in the semi-direct product $\mathbb{C}^*\ltimes\mathbb{C}$ is simply the
    multiplication action.)
 \end{prop}

  \begin{prop}{\label{P:nonsplit}}
     The 2-group $\PGL(d,d)$ is given by the following crossed-module:
         $$\xymatrix@=10pt@M=6pt{[\alpha \: \mathbb{C}^*
      \ar[rr]^{(\lambda^{d},\lambda^{d})}  &&
                                                \GL_2].   }$$
    We have $\pi_1\PGL(d,d)\cong \PGL_2$ and  $\pi_2\PGL(d,d)\cong \mathbb{Z}_d$.
  \end{prop}

 \section{Classification of spherical Deligne-Mumford curves}{\label{S:spherical}}

 Having described the
 structure of the 2-group of automorphisms of a weighted projective
 line $\PP(m,n)$, we can now classify spherical Deligne-Mumford curves.

 Let $\X$ be a spherical Deligne-Mumford curve.
 We know that the universal cover of $\X$ is a weighted projective
 line $\PP(m,n)$ (Proposition \ref{P:simpcntd2}). The deck transformations
 of the covering map $\PP(m,n) \to \X$ give us a homotopy class
 of map of 2-groups  $\pi_1\X \to \PGL(m,n)$ that is unique up to
 conjugation. Conversely, given such a conjugacy class, we can
 recover $\X$ up to isomorphism \cite{homotopy1}. So, what we
 have to do is to classify conjugacy classes of (homotopy classes of)
 maps $\pi_1\X \to \PGL(m,n)$.

Before presenting the classification result, let us explain how
this 2-group theoretic approach fits with the results of Section
\ref{S:UDMG}. We show how to derive the exact sequence of
Proposition \ref{P:spherical2} from the map $\pi_1\X \to
\PGL(m,n)$. (Note: the numbers $m$ and $n$ here are different
 from those of Proposition \ref{P:spherical2}.)

 By Proposition \ref{P:spherical2}, the universal cover of $\C$ is
$\PP(\frac{m}{d},\frac{n}{d})$, where $d=(m,n)$. Let
$\mathfrak{P}=\pi_1\X$, viewed as a 2-group, and let
$\mathfrak{H}$ be the  kernel of the map
 $\mathfrak{P}\to \PGL(m,n)$.
  After passing to classifying spaces and
 taking homotopy groups, we obtain the following fiber homotopy exact sequence:
 {\small
  $$\xymatrix@=7pt{1 \ar[r] & \pi_2\mathfrak{H} \ar[r] & \pi_2\mathfrak{P}
                       \ar[r]\ar@{=} [d] & \pi_2\PGL(m,n)  \ar@{=} [d]
       \ar [r] & \ar[r] \pi_1\mathfrak{H} \ar[r] & \pi_1\mathfrak{P} \ar@{=} [d]
                                       \ar[r] & \pi_1\PGL(m,n)  \\
            &&  1 &\Zd & & \pi_1\X & &}$$}
   It follows from the above exact
  sequence  that $\pi_2\mathfrak{H}$ is trivial.
  By Proposition \ref{P:automorphism}, we have
  $\pi_1\PGL(m,n)\cong\PGL(\frac{m}{d},\frac{n}{d})$. Via
  this isomorphism, the image of $\pi_1\X$ is identified with
  $\pi_1\C \subset \PGL(\frac{m}{d},\frac{n}{d})$. So the above exact
  sequence takes the following form:
    $$ 1 \to \Zd \to \pi_1\mathfrak{H} \to \pi_1\X \to \pi_1\C \to 1.$$
  It is easy to guess what $\pi_1\mathfrak{H}$ is; it is simply
  $H$, and the above exact sequence is isomorphic to
   that of Proposition \ref{P:spherical2}.

Now we get back to the classification of spherical Deligne-Mumford
curves. To do so, we
 combine Theorem \ref{T:maps} of the Appendix with the results
 of the preceding sections   to classify conjugacy
 classes of homotopy classes of maps $\Gamma \to \PGL(m,n)$ for
 any discrete group $\Gamma$.  As a byproduct of this
 classification, we  obtain an explicit
 (and canonical) description of a spherical Deligne-Mumford curve
 as a quotient of
 $\mathbb{C}^2-\{0\}$.

\vspace{0.1in}

 We consider two cases.

 \subsection{Universal cover $\PP(m,n)$, $m\neq n$}
 Let $d=\gcd(m,n)$.
  We saw in Section \ref{SS:explicit} that $\PGL(m,n)$ is split
  in this case (Proposition \ref{P:split}). We fix the splitting as in
  Section \ref{SS:explicit}.
  This immediately implies the following proposition.

 \begin{prop}{\label{P:mn}}
    For any group $\Gamma$,  there is a natural
    bijection
  $$ Hom_{\Ho(2Gp)}(\Gamma,\PGL(m,n))\cong
     \left\{ \begin{array}{ll}
       \Hom(\Gamma,\mathbb{C}^*)\times H^2(\Gamma,\mu_d)
                         & \text{if $m\nmid n$} \\
      \Hom(\Gamma,\mathbb{C}^*\ltimes\mathbb{C})\times H^2(\Gamma,\mu_m)
                         & \text{if $m\mid n$}
      \end{array}\right.$$
    where the the cohomology groups are with trivial coefficients and
    the semi-direct product is formed from the standard action of $\mathbb{C}^*$
    on $\mathbb{C}$.
 \end{prop}

 \begin{proof}
  Follows from Corollary \ref{C:split} of Appendix.
\end{proof}

To characterize Deligne-Mumford curves $\X$ with fundamental group
$\Gamma$ and universal cover ${\mathcal P}(m,n)$, what we need is
to classify conjugacy classes of elements in
$Hom_{\Ho(2Gp)}(\Gamma,\PGL(m,n))$. First we prove a lemma.

\begin{lem}{\label{L:conjugate}}
  For every
 group homomorphism $\psi\:\Gamma\to{\mathbb C}^\ast\ltimes\mathbb{C}$
 there is a unique character $\chi \: \Gamma \to \mathbb{C}^*$ that is
 conjugate to $\psi$. (Here we think of $\mathbb{C}^*$ as the first component
 of $\mathbb{C}^*\ltimes\mathbb{C}$.)
\end{lem}

\begin{proof}
 Every  group homomorphism $\psi\:\Gamma\to{\mathbb C}^\ast\ltimes\mathbb{C}$
 is of the form $(\chi,\eta)$, where $\chi$ is a
 character of $\Gamma$ and $\eta(u)=a(1-\chi(u))$ for some $a \in \mathbb{C}$.
 (If we think of ${\mathbb C}^\ast\ltimes\mathbb{C}$ as
 the group of affine transformations of $\mathbb{C}$, this is saying that
 the action of
 $\psi$  is  rotation
  around the point $a \in \mathbb{C}$, via the character
 $\chi$.) Using the element $(1,a) \in {\mathbb C}^\ast\ltimes\mathbb{C}$
 we can conjugate this with $\chi$ itself (in other words, we are sending
 $a$ to the origin).
\end{proof}

\begin{prop}{\label{P:conjugacy}}
   Assume $m\neq n$. Then the conjugacy classes of homotopy classes of
   maps from $\Gamma$ to $\PGL(m,n)$ are in bijection with
     $$\Hom(\Gamma,\mathbb{C}^*)\times H^2(\Gamma,\mu_d),$$
     where $d=\gcd(m,n)$.
\end{prop}

\begin{proof}
   If non of $m$ and $n$ divides the other, every two conjugate
  elements in $Hom_{\Ho(2Gp)}(\Gamma,\PGL(m,n))$ are in fact equal.

  When $m$ divides $n$, we know that
  $\pi_1\PGL(m,n)\cong {\mathbb C}^\ast\ltimes\mathbb{C}$.
  Lemma \ref{L:conjugate} implies that every conjugacy class in
  $$\Hom_{\Ho(\mathbf{2Gp})}(\Gamma,\PGL(m,n))\cong
       \Hom(\Gamma,\mathbb{C}^*\ltimes\mathbb{C})\times H^2(\Gamma,\mu_m)$$
  can be conjugated to an element $(\chi,c)$ in
   $$ \Hom(\Gamma,\mathbb{C}^*)\times H^2(\Gamma,\mu_m)\subseteq
      \Hom(\Gamma,\mathbb{C}^*\ltimes\mathbb{C})\times H^2(\Gamma,\mu_m).$$
  It also implies that $\chi$ is unique. So, we need to show that,
  if $(\chi,c)$ and $(\chi,c')$ are conjugate, then $c=c'$.
  In other words, we have to show that, for every $a \in G_{m,n}$
  such that $c_{\bar{a}}\circ\chi=\chi \: \Gamma \to \pi_1\mfG$, the
  induced conjugation action
    $$a_* \: \Hom_{\Ho(\mathbf{2Gp})}(\Gamma,\PGL(m,n))_{\chi}
     \to \Hom_{\Ho(\mathbf{2Gp})}(\Gamma,\PGL(m,n))_{\chi}$$
  is trivial. (Note that $a_*$ is not a priori a group homomorphisms.)

    By Corollary \ref{C:conjugation} of Appendix,
  it is enough to show that $a_*$ has a fixed point.  To obtain such a
  fixed point, consider the splitting $\pi_1\PGL(m,n) \to G_{m,n}$,
  and use it to lift $\chi$ to $\tilde{\chi} \: \Gamma \to G_{m,n}$
  We claim that $\tilde{\chi}$ is a fixed point of $a_*$.

  We may assume that $\chi$ is non-trivial.
  Notice that
  $\bar{a}\in \pi_1\PGL(m,n)\cong\mathbb{C}^*\ltimes\mathbb{C}$ commutes
  with every element in the image of $\chi$.
  It follows that $\bar{a} \in \mathbb{C}^*\subset\mathbb{C}^*\ltimes\mathbb{C}$.
  Therefore, $a$ is in the inverse image of $\mathbb{C}^* \in G_{m,n}$.
  This inverse image is isomorphic to $\mathbb{C}^*\times\mathbb{C}^*$,
  which is abelian.
  Since  the image of $\tilde{\chi}$ also lands in this abelian subgroup,
  conjugation by $a$ leaves $\tilde{\chi}$ invariant, which is what
  we wanted to prove.
\end{proof}

 By Proposition \ref{P:conjugacy}, every Deligne-Mumford curve $\X$ with
 universal cover $\PP(m,n)$ such that
 $\pi_1\X\cong\Gamma$ gives rise
 to an element in
 $\Hom(\Gamma,\mathbb{C}^*)\times H^2(\Gamma,\mu_d)$.
 Of course this element depends on the choice of the isomorphism
 $\pi_1\X\risom\Gamma$, so we have a bijection between
 the set of isomorphism classes
 of Deligne-Mumford curves whose fundamental groups is isomorphic
 to $\Gamma$ and the set
    $$\big(\Hom(\Gamma,\mathbb{C}^*)
        \times H^2(\Gamma,\mu_d)\big)/\Aut(\Gamma).$$

 Taking into account that elements in $ H^2(\Gamma,\mu_d)$ can
 be thought of as (isomorphisms classes of) central extensions
 of $\Gamma$ by $\mu_d$, this can be rephrased as follows.

\begin{thm}{\label{T:mn}}
  Let $m\neq n$ be given positive integer numbers, and let
  $d=\gcd(m,n)$. Then, the isomorphism classes of Deligne-Mumford
  curves $\X$ with universal cover $\PP(m,n)$ are in natural bijection with
  isomorphism
  classes of pairs $(K,\chi)$, where $K$ is a finite group containing
  $\mu_d$ as a central subgroup and
  $\chi \: K/\mu_d \to \mathbb{C}^*$ is a character.
\end{thm}

Another way to reformulate  Proposition \ref{P:conjugacy} is as
follows. Let $G_{m,n}' \subset G_{m,n}$ be the subgroup
$\mathbb{C}^*\times\mathbb{C}^*$. Form the  sub-crossed-module
$\PGL(m,n)'=[\varphi \: \mathbb{C}^* \to G_{m,n}']$. Proposition
\ref{P:conjugacy} is then equivalent to the following.

\begin{prop}{\label{P:conjugacy'}}
  Every conjugacy class in
   $\Hom_{\Ho(\mathbf{2Gp})}(\Gamma,\PGL(m,n))$ has a unique
   representative in $\Hom_{\Ho(\mathbf{2Gp})}(\Gamma,\PGL(m,n)')$.
   In other words, the set of conjugacy classes of elements in
   $\Hom_{\Ho(\mathbf{2Gp})}(\Gamma,\PGL(m,n))$ is in  natural
   bijection with $\Hom_{\Ho(\mathbf{2Gp})}(\Gamma,\PGL(m,n)')$.
\end{prop}

 \begin{cor}{\label{C:mn}}
 All Deligne-Mumford curves whose universal cover is ${\mathcal P}(m,n)$
 are given by triples $(\iota,E,\rho)$, where
  \begin{itemize}
   \item[$\mathbf{i.}$]
    $\iota:{\mathbb C}^\ast\hookrightarrow E$ is a
     central injection,
   \item[$\mathbf{ii.}$] $\rho:E\to \mathbb{C}^*\times\mathbb{C}^*$ is a
   group homomorphism,
  \end{itemize}
 such that the following two conditions are satisfied:

  \begin{itemize}
   \item[$\mathbf{a.}$]  $\rho\circ\iota:{\mathbb C}^\ast\to
     \mathbb{C}^*\times\mathbb{C}^*$
      is the diagonal $(m,n)$-th power map,
   \item[$\mathbf{b.}$] $E/{\mathbb C}^\ast$ is finite.
  \end{itemize}
 \noindent The curves given by $(\iota,E,\rho)$ and
 $(\iota',E',\rho')$ are isomorphic if and only if there is an isomorphism
 $f \: E \to E'$  such that
 $\rho'\circ f=\rho$.
\end{cor}

\begin{proof}In view of Proposition \ref{P:conjugacy'},
    this is just a special case of Theorem \ref{T:maps}
    of Appendix, with $\mfG=\PGL(m,n)'$.
\end{proof}

Let us describe how from a given  pair $(K,\chi)$ as in Theorem
\ref{T:mn}  we can reconstruct $\X$. What we will actually do is
to produce a triple $(\iota,E,\rho)$ as in  Corollary \ref{C:mn};
we will be implicitly using Corollary \ref{C:maps}  of Appendix.
 Let
$E:=\mathbb{C}^*\times_{\mu_d}K$ be the push out of $K$ along the
inclusion map $\mu_d \to \mathbb{C}^*$. Let
 $r$ and  $s$ be integers such that $sm+rn=d$ (as in Section \ref{SS:explicit}).
Extend the weighted action of $\mathbb{C}^*$ on $\mathbb{C}^2$ to
$E$
 by sending $(\lambda,u)\in E$ to the diagonal matrix
  $$\left( \begin{array}{cc} \chi(u)^{-r}\lambda^m  & 0 \\
                  0  &  \chi(u)^s\lambda^n \end{array}\right) $$
(This corresponds to the map $\rho$ of Corollary \ref{C:mn}.)
Then, $\X\cong[\mathbb{C}^2-\{0\}/E]$. This is an $H$-gerbe over
${\mathcal C}=[{\mathcal P}(\frac{m}{d},\frac{n}{d})/\chi(K)]$,
where $\chi(K) \subset \mathbb{C}^*$ is acting by rotations around
the origin, and $H=\ker\chi$.  The fundamental group of $\X$ is
isomorphic to $K/\mu_d$.

To get a better picture, let us complete the exact sequence of
Proposition \ref{P:spherical2}.$\mathbf{i}$ as follows:
   $$\xymatrix@C=10pt@R=8pt@M=5pt{ &&  K \ar[rd] &&&  \\
      1 \ar[r] & \mu_d \ar[r] &  H \ar[r] \ar@{^{(}->}[u]&
                \pi_1\X \ar[r] & \pi_1\C  \ar[r] & 1.      }$$

 \subsection{Universal cover $\PP(d,d)$}

  This case is more complicated than the previous one since
  $\PGL(d,d)$ is no longer split. As a consequence, the classification
  result is less explicit. The number $d$ is fixed throughout this section.

 \begin{prop}{\label{P:kk}}
    For any  discrete group $\Gamma$  there is a natural
    bijection between
    $Hom_{\Ho(2Gp)}(\Gamma,PGL(d,d))$ and the set of all
    isomorphism classes of diagrams
     $$\xymatrix@C=-4pt@R=10pt@M=10pt{
         & \mathbb{C}^*\ar[ld]\ar[rd]^{(\lambda^d,\lambda^d)}    &   \\
       E \ar[rr]_(0.45){\rho}                  &     &   \GL_2        }$$
    Here $E$ is a central extension of $\Gamma$ by ${\mathbb C}^*$,
    $\rho$ is a ${\mathbb C}$-representation of $E$ and the triangle
    commutes. (Isomorphism has the obvious meaning; see Theorem \ref{T:maps}.)
 \end{prop}

 \begin{proof}
    This follows from Theorem \ref{T:maps} of Appendix.
 \end{proof}

  Thus, every pair $(E,\rho)$ defines a Deligne-Mumford curve with
universal cover $\PP(d,d)$. Conversely, every such Deligne-Mumford
curve defines a pair $(E,\rho)$, that is unique up to conjugation
(this means, for any other $(E',\rho')$ there is an isomorphism $f
\: E \to E'$ and an element $a \in GL_2$ such that
 $\rho'\circ f=c_a\circ\rho$,
 where $c_a \: GL_2 \to GL_2$ is conjugation by $a$).

 Let us explain how this correspondence works.
Note that given $(E,\rho)$, taking $\pi_1$ gives a homomorphism of
groups $\Gamma\to PGL_2$.  Let $\overline{\Gamma}$ denote the
image. (Of course, $\overline{\Gamma}$ has to be `platonic'.)
Then, the Deligne-Mumford curve corresponding to $(E,\rho)$ is
$${\mathcal X}=[{\mathbb C}^2-\{0\}/E].$$
It is an $H$-gerbe over the orbifold ${\mathcal
C}=[\mathbb{P}^1/\overline{\Gamma}]$, where  $H=\ker\rho$.
Moreover, we have $\pi_1{\mathcal X}=\Gamma$ and
$\mathcal{X}=[\PP(d,d)/\Gamma]$.

The Deligne-Mumford curve, its orbifold and their universal covers
fit in the following commutative diagram:
   $$\xymatrix@=12pt@M=10pt{ \PP(d,d) \ar[r]^(0.52){B\mu_d}\ar[d]_{\Gamma}
                            &   {\mathbb P}^1 \ar[d]^{\overline{\Gamma}}  \\
           \X  \ar[r]_(0.38){BH}  &  [{\mathbb P}^1/{\overline{\Gamma}}]    }$$
(The arrows are labelled with their fibers.)

\begin{cor}{\label{C:kk}}
 All Deligne-Mumford curves whose universal cover is ${\mathcal P}(d,d)$
 are given by triples $(\iota,E,\rho)$, where
  \begin{itemize}
   \item[$\mathbf{i.}$]
    $\iota:{\mathbb C}^\ast\hookrightarrow E$ is a
     central injection,
   \item[$\mathbf{ii.}$] $\rho:E\to GL_2$ is a 2-dimensional representation,
  \end{itemize}
 such that the following two conditions are satisfied:

  \begin{itemize}
   \item[$\mathbf{a.}$]  $\rho\circ\iota:{\mathbb C}^\ast\to GL_2$
      is the diagonal $(d,d)$-th power map,
   \item[$\mathbf{b.}$] $E/{\mathbb C}^\ast$ is finite.
  \end{itemize}
 \noindent The curves given by $(\iota,E,\rho)$ and
 $(\iota',E',\rho')$ are isomorphic if and only if there is an isomorphism
 $f \: E \to E'$ and an element $a \in GL_2$ such that
 $\rho'\circ f=c_a\circ\rho$,
 where $c_a \: GL_2 \to GL_2$ is conjugation by $a$.
\end{cor}

The last part of the above corollary is simply saying that the
curves given by $(\iota,E,\rho)$ and $(\iota',E',\rho')$ are
isomorphic if and only if the corresponding maps in
$\Hom_{\Ho(\mathbf{2Gp})}(\Gamma,\PGL(d,d)$ are conjugate (see
page 38). This is unfortunately not as satisfactory as Corollary
9.6, because there, thanks to Proposition 9.3, we have a
recognition principle for conjugacy classes in in
$\Hom_{\Ho(\mathbf{2Gp})}(\Gamma,\PGL(m,n)$; in fact, Proposition
9.3 gives us a canonical representative for each conjugacy class.
We do not have such a clean picture in the case of  $\PGL(d,d)$,
making it more difficult to decide whether given two triples
$(\iota,E,\rho)$ and $(\iota',E',\rho')$ are conjugate or not. We
address this issue in more detail in the next subsection.

\subsection{Conjugacy classes in $\Hom_{\Ho(\mathbf{2Gp})}(\Gamma,\PGL(d,d))$}

We describe conjugacy classes of elements in
$\Hom_{\Ho(\mathbf{2Gp})}(\Gamma,\PGL(d,d))$. We do our analysis
by fixing a group homomorphism $\chi \: \Gamma \to \PGL_2$ and
studying the conjugacy classes in
 $\Hom_{\Ho(\mathbf{2Gp})}(\Gamma,\PGL(d,d))_{\chi}$.
We will assume throughout this section that this set is non-empty.

It  turns out that conjugacy classes have length 1,2 or 4, but
most of the time 1.
 We begin with an  easy proposition.

\begin{prop}{\label{P:centralizer}}
  Let $\chi \: \Gamma \to \PGL_2$ be a group homomorphisms, and
  let $I \subset PGL_2$ be its image. Let $C$ be the centralizer of
  $I$ in $\PGL_2$. Consider the action of $C$ on
  $\Hom_{\Ho(\mathbf{2Gp})}(\Gamma,\PGL(d,d))_{\chi}$ in which $a
  \in C$ acts by $c_{\tilde{a}}$, where $\tilde{a} \in  \GL_2$ is
  an  arbitrary lift of $a$. (Note that this is independent of the
  choice of the lift $\tilde{a}$.) Then, the orbits of this action
  are precisely the conjugacy classes of elements in
  $\Hom_{\Ho(\mathbf{2Gp})}(\Gamma,\PGL(d,d))_{\chi}$.
\end{prop}

\begin{proof}
  Obvious.
\end{proof}

To classify conjugacy classes of elements in
$\Hom_{\Ho(\mathbf{2Gp})}(\Gamma,\PGL(d,d))_{\chi}$ we consider
two cases: the image $I$ of $\chi$ is dihedral or not. The latter
case is easier.

\begin{prop}{\label{P:nonD}}
   Assume that the image $I\subset \PGL_2$ of $\chi$ is not
   dihedral. Then, every two conjugate elements in
   $\Hom_{\Ho(\mathbf{2Gp})}(\Gamma,\PGL(d,d))_{\chi}$ are
   actually equal.
\end{prop}

\begin{proof}
  If $I$ is not cyclic, a case by case inspection shows
    that  the centralizer of $I$ in $PGL_2$ is trivial.
    So the claim is true in these cases.
    In the case where $I$ is cyclic, we may assume that
    $I$ is a group of rotations around the origin. In this case,
    the centralizer of $I$ in $\PGL_2$ is the entire group $\mathbb{C}^*$
    of rotations    around the origin. The inverse image of this group in $\GL_2$
    is isomorphic to $\mathbb{C}^*\times\mathbb{C}^*$, which is abelian.
    So the claim is also true in this case.
\end{proof}

In the rest of this section we assume that $I$ is dihedral.

\begin{prop}{\label{P:dihedral}}
  Assume $I$ is dihedral, say $I\cong D_n=\mu_n\rtimes\mathbb{Z}_2$.
  Let $C$ be the centralizer of $I$ in $\PGL_2$. We have the
  following possibilities:
   \begin{itemize}
     \item  if $n$ is odd, then $C=\{1\}$.
    \item   if  $n$ is even,
        and $n>2$, then $C=\mu_2\subset \mu_n\subset D_n$.
    \item   if $I\cong D_2\cong\mathbb{Z}_2\times\mathbb{Z}_2$,
          then $C=I$.
  \end{itemize}
\end{prop}

\begin{proof}
  Easy. Left to the reader.
\end{proof}

We will now describe the action of $C$ on
$\Hom_{\Ho(\mathbf{2Gp})}(\Gamma,\PGL(d,d))_{\chi}$. We assume
$I=D_n$, $n$ even.
 For a non-trivial  element $a \in C$, define the
 map $a^* \: \Gamma \to \mathbb{Z}_2$ to be the composite $q\circ\chi$,
 where $q \: D_n \to \mathbb{Z}_2$ is the map obtained by killing $\mu_n
 \subset D_n$, when $n>2$. When $n=2$, $q$ is $D\to D/a=\mathbb{Z}_2$.

 \begin{lem}{\label{L:a}}
    Setup being as above, conjugation by $a\in C$ sends
    a triple
    $(\iota,E,\rho)$  in
    $\Hom_{\Ho(\mathbf{2Gp})}(\Gamma,\PGL(d,d))_{\chi}$
    to the triple $(\iota,E,\rho')$, where $\rho'=\rho  a^*$
    (pointwise product). Here, we have used the same notation
    $a^*$
    for the precomposition of $a^*$ with the projection
    map $E \to \Gamma$, and we have also identified $\mathbb{Z}_2$ with
    the corresponding scalar matrices.
 \end{lem}

\begin{proof}
    We may assume that $I=D_n=\mu_n\rtimes\mathbb{Z}_2$
     is the standard Dihedral group at the origin, with $n$ even.
     Here we have identified  $\mu_n$ with the group of matrices
       $$\left( \begin{array}{cc} \xi   & 0 \\
                  0  & 1 \end{array}\right)  \in \PGL_2 $$
    and the generator of $\mathbb{Z}_2$ is the matrix
        $$\left( \begin{array}{cc} 0   & 1 \\
                  1  &  0 \end{array}\right)   \in \PGL_2$$
    First we analyze   the conjugation action of $a$ on
    $J$, where  $J \in \GL_2$ is the inverse image of $I$. This is a
    central $\mathbb{C}^*$ extension of $I$. We only do the case
    $n>2$; the case $n=2$ is  similar.

     In this case  the
     centralizer of $I$ is $\mu_2$. We assume $a$ is non-trivial and
     choose the lift
     $\tilde{a}$ of $a$ to be the matrix
       $$\left( \begin{array}{cc} -1  & 0 \\
                  0  & 1 \end{array}\right) \in \GL_2$$
     It is easy to see that the conjugation action of this matrix
     leaves any matrix in $J$ that lies above $\mu_n \subset I$
     invariant, whereas it sends any other matrix $M$ in $J$ to
     $-M$; to see this, it is enough to take $M$ to be the following
     matrix:
        $$\left( \begin{array}{cc} 0  & 1 \\
                  1  & 0 \end{array}\right) $$
     Now, $\rho$ factors through $J$, so it follows from
     the above description that the conjugation action of $a$
      sends $\rho$ to $\rho a^*$.
     This proves the assertion.
\end{proof}

 We can think of $a^*$ as a character of $\Gamma$;
 that is, $a^* \in \Gamma^*$, where $\Gamma^*$ stands for the character group
 of $\Gamma$. This construction actually gives us an
 injection of $C$ in $\Gamma^*$ which sends $a$ to $a^*$. We
 denote the image by $C^* \subseteq \Gamma^*$. This is a subgroup
 of order 1,2 or 4.

\begin{prop}{\label{P:k}}
   In Lemma \ref{L:a}, $(\iota,E,\rho)$ and $(\iota,E,\rho')$
   are isomorphic, if and only if $a^*$ is in the image of the $d^{th}$
   power homomorphism $\Gamma^* \to \Gamma^*$.
\end{prop}

\begin{proof}
   First assume $a^*=\theta^d$, for some character $\theta \:
   \Gamma\to \mathbb{C}^{*}$. Define $f \: E \to E$ by
     $f(x)=x\iota(\theta(\bar{x}))$, where $\bar{x}$ denotes
     the image of $x$ in $\Gamma$ under the projection map
     $E \to \Gamma$.
     It is obvious that $f$ is an extension automorphism.
     It also satisfies $\rho'=\rho\circ f$, where
     $\rho':=c_a\circ\rho$ is the map described in Lemma
     \ref{L:a}. This implies that the triples
     $(\iota,E,\rho)$ and $(\iota,E,\rho')$ are isomorphic.

     Conversely, assume $(\iota,E,\rho)$ and $(\iota,E,\rho')$
     are isomorphic, and let the isomorphism be induced by some
     $f \: E \to E$. By definition, $f$ commutes
     with $\iota$ and the projection map to $\Gamma$,
     and it satisfies the relation $\rho'=\rho\circ f$.
     For an element $x \in E$ define $\theta(x)=x^{-1}f(x)$.
     This is an element in the image of $\iota$, so it gives us a
     map $\theta \: E \to \mathbb{C}^*$.
     In fact, since
     $\iota(\mathbb{C})$ is a central subgroup of $E$, it follows
     that $\theta$ is a group homomorphism.
     The relation $\rho'=\rho\circ f$ implies that $a^*=\theta^d$.
\end{proof}

\begin{cor}{\label{C:D}} Let $D=C^*/C^*\cap(\Gamma^*)^d$.
   Each orbit of $\Hom_{\Ho(\mathbf{2Gp})}(\Gamma,\PGL(d,d))_{\chi}$
   under conjugation action is a set of order $|D|$.
\end{cor}

\begin{proof}
  This follows immediately from Proposition \ref{P:k}.
\end{proof}

Remark that $|D|=1,2$ or 4. Except when $n$ and $d$ are both even,
we have $|D|=1$, and $|D|=4$ can only happen when $n=2$.

Summarizing the above discussion, if the image $I$ of the
homomorphism $\chi \: \Gamma \to \PGL_2$ is not a dihedral group,
then any two  elements in
$\Hom_{\Ho(\mathbf{2Gp})}(\Gamma,\PGL(d,d))_{\chi}$ that are
conjugate are indeed equal. If $I$ is dihedral, the conjugacy
classes have size 1, 2 or 4 (the first case always happens if at
least one of $d$ or $n$ is odd, the last case can only happen when
$I=D_2$). To determine which is the case, one has to compute the
characters $a^* \: \Gamma \to \mathbb{Z}_2$ as described above and
determine which ones are $d^{th}$ power. The ones that are
$d^{th}$ power induce trivial conjugation action on
$\Hom_{\Ho(\mathbf{2Gp})}(\Gamma,\PGL(d,d))_{\chi}$.

\section{Appendix: 2-groups}{\label{A:2gp}}

In this appendix we gather some results on 2-groups and
crossed-modules that are used throughout the paper. Our main
reference is \cite{2Gp}. More can be found in  \cite{Baez} (and
references therein) and \cite{Brown}.

\vspace{0.1in} \noindent{\bf Quick recall on 2-groups.}
 A {\bf 2-group} $\mathfrak{G}$ is a  group object in the category of groupoids.
 Alternatively, we can define a 2-group to be a groupoid object in the category
 of groups, or  also, as a (strict) 2-category with one object in
 which all  1-morphisms and 2-morphisms are invertible (in the strict sense).
 If we require that the 1-morphisms are only equivalences and not
 necessarily strictly invertible, what we obtain is called a
 {\em weak 2-group}.\footnote{In the case that we are interested in (namely,
 automorphisms of stacks)
 associativity
 holds strictly, so we have not  weakened associativity in our definition
 of a weak 2-group. That is why  we are sometimes sloppy and drop
 the adjective weak.} \footnote{Fortunately, we do not have to worry
 about the general strictification machinery for weak 2-groups, because
 the only weak 2-group we encounter indeed comes with a natural strictification
 (Theorem \ref{T:2-aut}).}

 \vspace{0.1in}
 \noindent{\em Example.} Self-equivalences of
 an object $\X$ in a 2-category (say a stack in the 2-category of stacks)
 naturally form a weak 2-group $\Aut(\X)$.

 \begin{rem}{\label{R:strict}}
  Weak 2-groups $\mathfrak{G}$ can be strictified in a functorial way
  (\cite{2Gp}), and strictification    retains
  all the (homotopical) properties of the original weak 2-group.
  For this reason in practice we will only need to work with strict
  2-groups.
 \end{rem}

 A {\em morphism}
 $f \: \mathfrak{G} \to \mathfrak{H}$ of 2-groups is a map of groupoids that
 respects the group operation on the nose.
 If we view $\mathfrak{G}$ and  $\mathfrak{H}$
 as 2-categories with one object, such $f$ is   nothing but  a strict 2-functor.
 Same definition applies to weak 2-groups (because we have not
 weakened associativity).

 To a 2-group $\mathfrak{G}$ we associate the groups
 $\pi_1\mathfrak{G}$ and $\pi_2\mathfrak{G}$ as follows. The group
 $\pi_1\mathfrak{G}$
 is the set of isomorphism classes of object of the groupoid $\mathfrak{G}$.
 The group structure on $\pi_1\mathfrak{G}$ is induces from the group
 operation of $\mathfrak{G}$. The group
 $\pi_2\mathfrak{G}$ is the group of automorphisms  of the identity
 object $e \in \mathfrak{G}$. This is an abelian group.
 Any  morphism
 $f \: \mathfrak{G} \to \mathfrak{H}$ of (weak) 2-groups induces homomorphisms
 on $\pi_1$ and $\pi_2$. We say such $f$ is
  an {\em equivalence}
 if both these maps are isomorphisms.
 ({\em Warning}: an equivalence
 need not have an inverse. Also, two equivalent
  2-groups may not be related by an equivalence, but just a zig-zag of
  equivalences.)

 We are usually interested
 in 2-groups up to equivalence, so we will think of a 2-group as an object
 in the {\em homotopy category} of 2-groups. This category is defined by taking
 the category of 2-groups   and
 inverting the  equivalences.
 For the reader's information, we point out that,
 there is a model structure on the category
 of 2-groupoids (see for
 instance \cite{Moerdijk}); 2-groups are the pointed connected objects in
 this category.

 \vspace{0.1in}
 \noindent{\bf Quick recall on crossed-modules.} A {\bf crossed-module}
 $\mathfrak{G}=[\varphi \: G_2 \to G_1]$
  is a pair of groups $G_1,G_2$, a group homomorphism $\varphi \: G_2 \to G_1$,
  and a (right) action of $G_1$ on $G_2$, denoted $-^g$,
  that lifts the conjugation action
  of $G_1$ on
  the image of $\varphi$ and descends the conjugation action of $G_2$ on itself.
  The kernel of $\varphi$ is a central (in particular abelian) subgroup of $G_2$
  and is denoted by $\pi_2\mathfrak{G}$. The image of $\varphi$ is a normal
  subgroup of $G_1$ whose cokernel is denoted by  $\pi_1\mathfrak{G}$.
  A  (strict) morphism  of crossed-modules is a pair of group homomorphisms
  which commute
  with the $\varphi$ maps and respect the action of $G_1$ on $G_2$.
  A  morphism is called an {\em equivalence} if
  it induces an isomorphism on $\pi_1$ and $\pi_2$.

\vspace{0.1in} \noindent{\bf Equivalence of 2-groups and
crossed-modules}
   There is a natural pair of inverse equivalences of categories between
   the category $\mathbf{2Gp}$ of 2-groups  and the category
   $\mathbf{CrossedMod}$ of crossed-modules.
   Furthermore,
   these functors preserve $\pi_1$ and $\pi_2$. Here is how these
   functors are defined.

   \vspace{0.1in}
   \noindent{\em Functor from 2-groups to crossed-modules.}
   Let $\mathfrak{G}$ be a 2-group. Let $G_1$ be the group of  objects of
   $\mathfrak{G}$, and let $G_2$ be the set of all arrows emanating from the
   identity object $e$; $G_2$ is also a group (namely, it is a subgroup
   of the group  of arrows of $\mathfrak{G}$).

   Define $\varphi \: G_2 \to G_1$ by
   $\varphi(\alpha):=t(\alpha)$,
   $\alpha \in G_2$.

   The action of $G_1$ on $G_2$ is given by conjugation. That is,
   given $\alpha \in G_2$ and $g \in G_1$, the action is given by
   $g^{-1}\alpha g$.

   \vspace{0.1in}
   \noindent{\em Functor from crossed-modules to 2-groups.}
   Let  $[\varphi \: G_2 \to G_1]$ be a crossed-module. Consider the
   groupoid $\mathfrak{G}$ whose underlying set of objects is $G_1$ and whose
   set of arrows is $G_1\times G_2$. The source and target maps are given by
   $s(g,\alpha)=g$, $t(g,\alpha)=g\varphi(\alpha)$.

   Two arrows $(g,\alpha)$ and $(h,\beta)$ such that $g\varphi(\alpha)=h$
   are composed to $(g,\alpha\beta)$.

   The group operation on the set of objects $\Ob\mathfrak{G}=G_1$ is naturally
   extended to a group operation on $\mathfrak{G}$ by setting
   $(g,\alpha)(h,\beta)=(gh,(\alpha^h)\beta)$, where $-^h$ stands for the
   action of $G_1$ on $G_2$.

  \begin{rem}{\label{R:synonym}}
     With the above equivalence of categories in mind, we use the
     words 2-group and crossed-module synonymously throughout the
     paper. Every discrete group $\Gamma$ gives rise to a 2-group in the
     obvious way. The crossed-module presentation for this group
     is $[1 \to \Gamma]$. This gives a full embedding of the
     category of groups in the category of 2-groups (or
     crossed-modules).
  \end{rem}

\subsection{Morphisms into 2-groups}

Let $\Gamma$ be a discrete group and $\mathfrak{G}$ a 2-group. We
are interested in determining the homotopy classes of morphisms
from $\Gamma$ to $\mathfrak{G}$. This arises in the study of group
actions on stacks. Namely, let $\X$ be a stack, and
$\mathfrak{G}=\Aut(\X)$ be the automorphism 2-group of $\X$
(rather, its strictification). Then the actions of $\Gamma$ on
$\X$ (up to a reasonable notation of equivalence) are in a natural
bijection with homotopy classes of maps from $\Gamma$ to
$\mathfrak{G}$.

\vspace{0.1in} \noindent{\em Caveat.} Strictly speaking, our use
of the expression `homotopy classes of morphisms from $\Gamma$ to
$\mathfrak{G}$' is not quite correct. What we really mean is
`elements in $\Hom_{\Ho(\mathbf{2Gp})}(\Gamma,\mathfrak{G})$'. The
point is that, not every element in
$\Hom_{\Ho(\mathbf{2Gp})}(\Gamma,\mathfrak{G})$ can be realized by
a strict morphism of  2-groups $\Gamma \to \mathfrak{G}$, but only
a {\em weak} one.

\vspace{0.1in} The following result is proved in \cite{2Gp}.

\begin{thm}[\oldcite{2Gp}, Corollary 6.9]{\label{T:maps}} Let $\mfG$ be a
  2-group and  let  $[\varphi \: G_2 \to G_1]$  be the
  crossed-module presentation
for it. Let $\Gamma$ be a discrete
  group. Then, the set $\Hom_{\Ho(\mathbf{2Gp})}(\Gamma,\mathfrak{G})$
  of homotopy classes of morphisms from $\Gamma$ to $\mfG$
  is in a natural bijection with the set of isomorphism classes of
  commutative diagrams
  of the form
   $$\xymatrix@C=-4pt@R=10pt@M=10pt{
         & G_2\ar[ld]\ar[rd]^{\varphi}    &   \\
       E \ar[rr]_(0.45){\rho}                  &     &   G_1         }$$
  where $E$ is an extension of $\Gamma$ by $G_2$, and for every $x \in E$
  and $g \in G_2$ we require that $g^{\rho(x)}=x^{-1}gx$. An
  isomorphism between two such diagrams $(E,\rho)$ and $(E',\rho')$ is,
  by definition, an isomorphism of extensions $f \: E' \to E$ (so
  it induces identity on $G_2$ and $\Gamma$) such that
  $\rho'=\rho\circ f$.
\end{thm}

In terms of actions of discrete groups on stacks, this can be
interpreted as follows. Let $\X$ be a stack (or, more generally,
an object in a 2-category). Take  $\mathfrak{G}$ to be a strict
model for the (weak) 2-group of automorphisms of $\X$, with
$[\varphi \: G_2 \to G_1]$ the crossed-module presentation for it.
The elements of $G_1$ are then to be thought of as a complete
collection of (2-isomorphism classes) of automorphisms of $\X$,
while the elements of $G_2$ parameterize 2-isomorphisms between
them. Of course, too few actions of $\Gamma$ on $\X$ can
 be actually obtained as group homomorphisms $\Gamma \to G_1$, as we
have to allow certain ``laxness''. What the theorem is saying is
that, there is, however, a canonical way of weaving the elements
of $G_2$ into $\Gamma$ to form a larger group $E$ so that  our
action now can be realized as an actual group homomorphism $E \to
G_1$, with $G_2 \subseteq E$ keeping track of the 2-isomorphisms.

For computational purposes, it is useful to have the cohomological
versions of Theorem \ref{T:maps} handy.

\begin{thm}[\oldcite{2Gp}, Proposition 7.6]{\label{T:cohomological}}
     Let $\mathfrak{G}$ be a 2-group, and let $\Gamma$
     be a discrete group. Fix a homomorphism
     $\chi \: \Gamma \to \pi_1\mathfrak{G}$. Consider the set
     $\Hom_{\Ho(\mathbf{2Gp})}(\Gamma,\mathfrak{G})_{\chi}$ of homotopy
     classes of maps $\Gamma \to \mfG$ inducing $\chi$ on $\pi_1$. Then,
     either this set is empty, or it admits a natural
     transitive action of
     $H^2(\Gamma,\pi_2\mathfrak{G})$.
     (Here, $\pi_2\mathfrak{G}$ is made into a $\Gamma$-module via $\chi$.)
\end{thm}

\begin{defn}{\label{D:split}}
   We say a 2-group $\mathfrak{G}=[G_2 \to G_1]$ is {\em split}
   if the projection map $G_1 \to \pi_1\mathfrak{G}$ admits a section.
   (Note: this definition is slightly different from the one
   used in \oldcite{2Gp}.)
\end{defn}

\begin{cor}{\label{C:split}}
  In Theorem \ref{T:maps}, assume the 2-group $\mfG$ is split, and fix a section
  $\sigma \: \pi_1\mfG \to G_1$. Then, there is a
  natural bijection
    $$\Hom_{\Ho(\mathbf{2Gp})}(\Gamma,\mathfrak{G})_{\chi}
                     \longleftrightarrow H^2(\Gamma,\pi_2\mathfrak{G}).$$
\end{cor}

In Theorem \ref{T:maps}, the  action of
$H^2(\Gamma,\pi_2\mathfrak{G})$ can be made completely explicit
(see \cite{2Gp}), and this allows us to  write down every homotopy
class $\Gamma \to \mfG$, provided that we are given one of them to
begin with. We will show how to do this in the special case where
$\mfG$ is split. First we need to introduce a notation.

 \begin{defn}{\label{D:crossed}}
Let $H$, $G$ and $K$ be  groups, each equipped with a right action
of $K$, the one on $K$ itself being conjugation. We denote the
actions by $-^k$. Assume we are given a $K$-equivariant diagram:
  $$\xymatrix@=12pt@M=10pt{
     H \ar[r]^p  \ar[d]_l  &   G   \\
          K            &       }$$
satisfying the  compatibility condition $g^{l(h)}=p(h)^{-1}gp(h)$,
for
 $h\in H$ and $g\in G$.
    The {\em semi-direct product} $K\ltimes_H G$ of $K$ and $G$ along $H$
    is defined to be $K\ltimes G/N$, where
      $$N=\big\{\big(l(h)^{-1},p(h)\big), \ h \in H \big\}.$$
 \end{defn}

\begin{cor}[\oldcite{2Gp}, Corollary 7.9]{\label{C:maps}}
  Suppose the natural projection $G_1 \to \pi_1\mfG$ has a
  section,
  and fix such a section $\sigma \: \pi_1\mfG \to G_1$.
   Then, every homotopy class $\Gamma \to \mfG$
  is uniquely characterized by  a pair $(K,\chi)$, where
  $\chi \: \Gamma \to \pi_1\mfG$
  is a group homomorphism and $K$ is (the isomorphism class of)
  an extension of $\Gamma$ with $\pi_2\mathfrak{G}$:
     $$1 \to \pi_2\mathfrak{G} \to K \to \Gamma  \to 1.$$
  Given such an extension we obtain  a diagram as in
  Theorem \ref{T:maps} by setting
  $E:=K \ltimes_{\pi_2\mathfrak{G}} G_2$
  and  $\rho(k,\alpha):=\sigma\big(\chi(\bar{k})\big)\varphi(\alpha)$,
  $(k,\alpha) \in E$.
  Here, the action of $K$ on $G_2$
  is obtained via $\sigma\circ\chi$ from that of $G_1$ on $G_2$.
  The notation $\bar{k}$ stands for the image of $k \in K$ in $\Gamma$.
\end{cor}

Again, we can interpret  this result in terms of group action on
stack as follows. First of all, to fix a section $\sigma \:
\pi_1\mfG \to G_1$ means to fix a  strict action of $\pi_1\mfG$ on
$\X$. The corollary now says that, to give an action of $\Gamma$
on $\X$ is the same as letting an extension $K$ of $\Gamma$ by
$\pi_2\mfG$ act strictly on $\X$, with $\pi_2\mfG \subseteq K$
acting trivially.

\subsection{Conjugations in 2-groups and crossed-modules} {\label{SS:conjugation}}

  Given a 2-group $\mathfrak{G}$ and an element $a$ in $G_1$ (the ``group of
  objects'') we  define the automorphism, $c_a\: \mfG \to \mfG$,
  called {\em conjugation} by $a$, to be the map that
  sends an object $h$ (respectively, and arrow $\alpha$)  to
  $a^{-1}ha$ (respectively, $a^{-1}\alpha a$). This is obviously
  a map of 2-groups. On the crossed-module $[\varphi \: G_2 \to G_1]$
  associated to $\mathfrak{G}$, the conjugation by $a$ acts as
  conjugation  by $a$ on $G_1$ and as $-^a$ on $G_2$.

  Using this we can define the notion of conjugacy classes of maps
  between 2-groups. This notion passes to the homotopy
  category.
  In particular, we can also talk about conjugacy classes
  of maps between  two objects in the homotopy category of
  2-groups.

  Let $\Gamma$ be a group, and assume we are given an element in
  $\Hom_{\Ho(\mathbf{2Gp})}(\Gamma,\mathfrak{G})$ represented by a
  diagram
  $$\xymatrix@C=-4pt@R=10pt@M=10pt{
         & G_2\ar[ld]_{\eta} \ar[rd]^{\varphi}    &   \\
       E \ar[rr]_(0.45){\rho}                  &     &   G_1         }$$
 as in Theorem \ref{T:maps}.
 For $a \in G_1$, the conjugation action
  $c_a \: \Hom_{\Ho(\mathbf{2Gp})}(\Gamma,\mathfrak{G}) \to
    \Hom_{\Ho(\mathbf{2Gp})}(\Gamma,\mathfrak{G})$
    sends this diagram to the diagram
    $$\xymatrix@C=-4pt@R=10pt@M=10pt{
         & G_2\ar[ld]_{\eta\circ c_a^{-1} } \ar[rd]^{\varphi}    &   \\
       E \ar[rr]_(0.45){c_a\circ\rho}                  &     &   G_1        }$$
 Here we have abused the notation and
 denoted the conjugation action of $a$ on $G_1$, as
 well as the action of $a$ on $G_2$, by $c_a$.

    Let us  examine the effect of the conjugation action from the
    cohomological point of view. For a given $\chi \: \Gamma \to \pi_1\mfG$,
    and an element $a \in G_1$,
    the conjugation $c_a \: \mfG \to \mfG$ induces a natural
    map
      $$a_* \: \Hom_{\Ho(\mathbf{2Gp})}(\Gamma,\mathfrak{G})_{\chi} \to
      \Hom_{\Ho(\mathbf{2Gp})}(\Gamma,\mathfrak{G})_{c_{\bar{a}}\circ\chi}.$$
    This map is equivariant in the following sense. Recall
     that
    $\Hom_{\Ho(\mathbf{2Gp})}(\Gamma,\mathfrak{G})_{\chi}$ is an
    $H^2(\Gamma,\pi_2\mfG)$-set (Theorem \ref{T:cohomological}). Similarly,
    $\Hom_{\Ho(\mathbf{2Gp})}(\Gamma,\mathfrak{G})_{c_{\bar{a}}\circ\chi}$ is
    an
    $H^2(\Gamma,\pi_2\mfG)_a$-set. The index $_a$ here is to emphasis
    that,
    in the latter cohomology group the $\Gamma$-action on $\pi_2\mfG$
    is different from the  one in the former cohomology group; the two actions are
    conjugate by the automorphism $-^a \: \pi_2\mfG \to \pi_2\mfG$.
    The automorphisms $-^a \: \pi_2\mfG \to \pi_2\mfG$,
    which is  $\Gamma$-equivariant for the two actions, induces an isomorphism
    $\alpha \: H^2(\Gamma,\pi_2\mfG) \to H^2(\Gamma,\pi_2\mfG)_a$.
    The map $a_*$ is an $\alpha$-equivariant map.

    \begin{cor}{\label{C:conjugation}}
      Notation being as above, assume $a$ has the following
      properties: 1) the action of $a$  on $\pi_2\mfG$ is trivial, 2)
      $c_a\circ\chi=\chi$ and
      there is a class in
      $\Hom_{\Ho(\mathbf{2Gp})}(\Gamma,\mathfrak{G})_{\chi}$ that is invariant
      under $a_*$. Then $a_*$ is identity.
    \end{cor}

    \begin{proof}
       We may assume $\Hom_{\Ho(\mathbf{2Gp})}(\Gamma,\mathfrak{G})_{\chi}$
       is non-empty.
       By hypothesis, the isomorphism
       $\alpha \: H^2(\Gamma,\pi_2\mfG) \to H^2(\Gamma,\pi_2\mfG)_a$ is
       identity. Since $a_*$ is an $\alpha$-equivariant self-map of a transitive
       $H^2(\Gamma,\pi_2\mfG)$-set, and since it
       has a fixed point, it must be identity.
    \end{proof}

\bibliographystyle{amsplain}
\bibliography{uniformization}
\end{document}

%% file: fig1.pstex_t
\begin{picture}(0,0)%
\includegraphics{fig1.pstex}%
\end{picture}%
\setlength{\unitlength}{3947sp}%
\begingroup\makeatletter\ifx\SetFigFont\undefined%
\gdef\SetFigFont#1#2#3#4#5{%
  \reset@font\fontsize{#1}{#2pt}%
  \fontfamily{#3}\fontseries{#4}\fontshape{#5}%
  \selectfont}%
\fi\endgroup%
\begin{picture}(5457,1853)(248,-1315)
\put(2476,-811){\makebox(0,0)[lb]{\smash{\SetFigFont{8}{9.6}{\rmdefault}{\mddefault}{\updefault}$G_2$}}}
\put(826,-736){\makebox(0,0)[lb]{\smash{\SetFigFont{8}{9.6}{\rmdefault}{\mddefault}{\updefault} $G_k$}}}
\put(2271,-106){\makebox(0,0)[lb]{\smash{\SetFigFont{8}{9.6}{\rmdefault}{\mddefault}{\updefault}$G_1$}}}
\put(4697,-182){\makebox(0,0)[lb]{\smash{\SetFigFont{8}{9.6}{\rmdefault}{\mddefault}{\updefault}$G_2$}}}
\put(4514,318){\makebox(0,0)[lb]{\smash{\SetFigFont{8}{9.6}{\rmdefault}{\mddefault}{\updefault}$G_1$}}}
\put(5446,-588){\makebox(0,0)[lb]{\smash{\SetFigFont{8}{9.6}{\rmdefault}{\mddefault}{\updefault}$G_k$}}}
\put(901,-136){\makebox(0,0)[lb]{\smash{\SetFigFont{8}{9.6}{\rmdefault}{\mddefault}{\updefault}$H$}}}
\put(4876,-1111){\makebox(0,0)[lb]{\smash{\SetFigFont{8}{9.6}{\rmdefault}{\mddefault}{\updefault}$H$}}}
\end{picture}

%% file: fig2.pstex_t
\begin{picture}(0,0)%
\includegraphics{fig2.pstex}%
\end{picture}%
\setlength{\unitlength}{4144sp}%
\begingroup\makeatletter\ifx\SetFigFont\undefined%
\gdef\SetFigFont#1#2#3#4#5{%
  \reset@font\fontsize{#1}{#2pt}%
  \fontfamily{#3}\fontseries{#4}\fontshape{#5}%
  \selectfont}%
\fi\endgroup%
\begin{picture}(2090,2050)(1108,-2670)
\put(1138,-1521){\makebox(0,0)[lb]{\smash{\SetFigFont{8}{9.6}{\rmdefault}{\mddefault}{\updefault}$G_1$}}}
\put(1108,-1804){\makebox(0,0)[lb]{\smash{\SetFigFont{8}{9.6}{\rmdefault}{\mddefault}{\updefault}$G_2$}}}
\put(1920,-2639){\makebox(0,0)[lb]{\smash{\SetFigFont{8}{9.6}{\rmdefault}{\mddefault}{\updefault}$G_k$}}}
\put(1623,-1786){\makebox(0,0)[lb]{\smash{\SetFigFont{8}{9.6}{\rmdefault}{\mddefault}{\updefault}$H$}}}
\put(1635,-1654){\makebox(0,0)[lb]{\smash{\SetFigFont{8}{9.6}{\rmdefault}{\mddefault}{\updefault}$H$}}}
\put(2136,-956){\makebox(0,0)[lb]{\smash{\SetFigFont{8}{9.6}{\rmdefault}{\mddefault}{\updefault}$H$}}}
\put(2430,-1222){\makebox(0,0)[lb]{\smash{\SetFigFont{8}{9.6}{\rmdefault}{\mddefault}{\updefault}$H$}}}
\put(2655,-2160){\makebox(0,0)[lb]{\smash{\SetFigFont{8}{9.6}{\rmdefault}{\mddefault}{\updefault}$H$}}}
\put(1956,-2188){\makebox(0,0)[lb]{\smash{\SetFigFont{8}{9.6}{\rmdefault}{\mddefault}{\updefault}$H$}}}
\end{picture}